\documentclass[11pt]{amsart}
\usepackage{amsxtra}
\theoremstyle{plain}

\newcommand{\cleqn}{\setcounter{equation}{0}}
\newcommand{\clth}{\setcounter{theorem}{0}}
\newcommand {\sectionnew}[1]{\section{#1}\cleqn\clth}

\newtheorem{theorem}{Theorem}[section]
\newtheorem{lemma}[theorem]{Lemma}
\newtheorem{definition-theorem}[theorem]{Definition-Theorem}
\newtheorem{proposition}[theorem]{Proposition}
\newtheorem{corollary}[theorem]{Corollary}
\newtheorem{definition}[theorem]{Definition}
\newtheorem{example}[theorem]{Example}
\newtheorem{remark}[theorem]{Remark}
\newtheorem{notation}[theorem]{Notation}
\newcommand \bth[1] { \begin{theorem}\label{t#1} }
\newcommand \ble[1] { \begin{lemma}\label{l#1} }

\newcommand \bpr[1] { \begin{proposition}\label{p#1} }
\newcommand \bco[1] { \begin{corollary}\label{c#1} }
\newcommand \bde[1] { \begin{definition}\label{d#1}\rm }
\newcommand \bex[1] { \begin{example}\label{e#1}\rm }
\newcommand \bre[1] { \begin{remark}\label{r#1}\rm }

\newcommand \bnota[1] {\begin{notation}\label{n#1}\rm }
\newcommand {\eth} { \end{theorem} }
\newcommand {\ele} { \end{lemma} }

\newcommand {\epr} { \end{proposition} }
\newcommand {\eco} { \end{corollary} }
\newcommand {\ede} { \end{definition} }
\newcommand {\eex} { \end{example} }
\newcommand {\ere} { \end{remark} }
\newcommand {\enota} { \end{notation} }

\newcommand \thref[1]{Theorem \ref{t#1}}
\newcommand \leref[1]{Lemma \ref{l#1}}
\newcommand \prref[1]{Proposition \ref{p#1}}
\newcommand \coref[1]{Corollary \ref{c#1}}

\newcommand \reref[1]{Remark \ref{r#1}}
\newcommand \lb[1]{\label{#1}}

\def \Cset {{\mathbb C}}

\def \O  {{\mathcal{O}}} 
\def \L  {{\mathcal{L}}}
\def \Ld  {{\mathcal{L}}(\d)}
\def \calS {{\mathcal{S}}}
\def \De {\Delta}   

\def \al {\alpha}

\def \la {\lambda}

\def \ra  {\rightarrow}           
           
\def \lra {\longrightarrow}

\def \la {\langle}
\def \ra {\rangle}

\def \Ad { {\mathrm{Ad}} }

\def \Lie { {\mathrm{Lie}} \,}

\def \g  {\mathfrak{g}}   
\def \h  {\mathfrak{h}}
\def \f  {\mathfrak{f}}
\def \n  {\mathfrak{n}}
\def \m  {\mathfrak{m}}
\def \b  {\mathfrak{b}}
\def \p  {\mathfrak{p}}
\def \r  {\mathfrak{r}}
\def \z  {\mathfrak{z}}

\def \u  {\mathfrak{u}}
\def \up  {\mathfrak{u}^\prime}

\def \q  {\mathfrak{q}}
\def \d  {\mathfrak{d}}

\def \l  {\mathfrak{l}}
\def \Up {U^\prime}
\def \a {\mathfrak{a}}

\DeclareMathOperator \Span { {\mathrm{Span}} }

\DeclareMathOperator \ad { {\mathrm{ad}} }

\DeclareMathOperator \Rank { {\mathrm{Rank}} }

\newcommand \Gr { {\mathrm{Gr}} }
\renewcommand{\qed}{\begin{flushright} {\bf Q.E.D.}\ \ \ \ \
                  \end{flushright} }    
\def \hs {\hspace{.2in}}
\def \lara {\la \, \, , \, \ra}
\def \ua {\underline{a}}

\def \ud {\underline{d}}
\def \ue {\underline{e}}
\def \Lgg {\L(\g \oplus \g)}
\def \gog {\g \oplus \g}
\def \Pill {\Pi_{\l_1, \l_2}}
\def \hoh {\h \oplus \h}
\def \lara {\la  \, , \, \ra}
\begin{document}
\setlength{\baselineskip}{1.2\baselineskip}
\title[Regular partitions of
Poisson manifolds]
{Group orbits and regular partitions of
Poisson manifolds}
\author[Jiang-Hua Lu]{Jiang-Hua Lu}
\address{
Department of Mathematics \\
Hong Kong University \\
Pokfulam Rd., Hong Kong}
\email{jhlu@maths.hku.hk}
\author[Milen Yakimov]{Milen Yakimov}
\address{
Department of Mathematics \\
University of California \\
Santa Barbara, CA 93106, U.S.A.}
\email{yakimov@math.ucsb.edu}
\date{}
\begin{abstract} We study a large class of Poisson manifolds,
derived from Manin triples, for which we construct explicit 
partitions into regular Poisson submanifolds by
intersecting certain group orbits.
Examples include all varieties ${\mathcal L}$ of Lagrangian subalgebras of
reductive quadratic Lie algebras $\d$ with Poisson structures defined by 
Lagrangian splittings of $\d$. In the special case of 
$\g \oplus \g$, where $\g$ is a complex semi-simple Lie algebra, we 
explicitly compute
the ranks of the Poisson structures on ${\mathcal L}$ defined by 
arbitrary Lagrangian splittings of ${\mathfrak g} \oplus {\mathfrak g}$. Such 
Lagrangian splittings have been classified by P. Delorme, and they contain 
the Belavin--Drinfeld splittings as special cases. 
\end{abstract}
\maketitle
\sectionnew{Introduction}\lb{intro}

Lie theory provides a rich class of examples of Poisson manifolds/varieties. 
In this paper, we study a  class of Poisson manifolds of the form
$(D/Q, \Pi_{\u, \up})$, where 
$D$ is an even dimensional connected real or complex Lie group whose Lie 
algebra
$\d$ is quadratic, i.e. $\d$ is equipped 
with a nondegenerate invariant symmetric bilinear form $\lara$; the closed 
subgroup
$Q$ of $D$ corresponds to a subalgebra $\q$ of $\d$ that
is coisotropic with respect to $\lara$, and $(\u, \up)$ is
a pair of complementary Lagrangian subalgebras of $\d$.
A Lie subalgebra $\l$ of $\d$ will be called Lagrangian 
if $\l^\perp = \l$ with respect to $\lara$. We will call such a 
splitting $\d = \u + \up$ a Lagrangian splitting.
The Poisson  
structure $\Pi_{\u, \up}$ is obtained from the $r$-matrix
\[
r_{\u, \up} = \frac{1}{2} \sum_{i=1}^n \xi_i \wedge x_i \in \wedge^2 \d,
\]
where $\{x_1, \ldots, x_n\}$
and $\{\xi_1, \ldots, \xi_n\}$ are pairs of dual bases of $\u$ and 
$\up$ with respect to $\lara$. We refer the reader to $\S$ \ref{sec_DQ}
for the precise definition of $\Pi_{\u, \up}$.

Let $U$ and $U^\prime$ be the connected 
subgroups of $D$ with Lie algebras $\u$ and $\up$
respectively. Our first main result, 
see \thref{qqqperp} and \prref{N-regular}, 
is that when
\begin{equation}
 [\q, \q] \subset \q^\perp,
\label{cond}
\end{equation}
all intersections of $U$ and $U^\prime$-orbits in $D/Q$ are 
regular Poisson submanifolds. In fact, if $N(\u)$ and $N(\up)$
denote the normalizers of
$\u$ and $\up$ in $D$ respectively, we also show that all
intersections of $N(\u)$ and $N(\up)$-orbits in $D/Q$ are 
regular Poisson submanifolds.
Note that the condition
\eqref{cond} is an intrinsic property of the coisotropic 
subalgebra $\q$ of $(\d, \lara)$ and does not depend on the 
Lagrangian splitting $\d = \u + \up$. Once this condition is verified
for a given $\q$, the above result provides ``regular'' partitions 
for the Poisson structures $\Pi_{\u, \up}$ on $D/Q$ for any 
Lagrangian splitting $\d = \u + \up$. 

Our second main result shows that the condition \eqref{cond} is satisfied 
when 
$\d$ is reductive and $\q$ is the normalizer subalgebra in $\d$ of any
Lagrangian subalgebra of $(\d, \lara)$. In fact, we show in
\prref{cond} that in this case 
$[\q, \q] = \q^\perp$. Let
$\L(\d, \lara)$ be the variety of Lagrangian subalgebras of 
$(\d, \lara)$. All $D$-orbits in $\L(\d, \lara)$
are of the form $D/N(\l)$, where $N(\l)$ is the normalizer 
subgroup in $D$ of an $\l \in \L(\d, \lara)$, and every Lagrangian
splitting $\d = \u + \up$ defines a Poisson structure $\Pi_{\u, \up}$ on
$\L(\d, \lara)$. As a corollary of the second main result we obtain that,
if $\d$ is an even dimensional reductive quadratic Lie algebra, then 
every non-empty intersection of an $N(\u)$-orbit and an $N(\up)$-orbit 
on $\L(\d, \lara)$ is a regular Poisson submanifold with respect to 
$\Pi_{\u, \up}$

In $\S$\ref{Lgg}, we take the special case when $\d = \gog$ for a 
complex semi-simple Lie algebra $\g$ and 
\[
\la (x_1, x_2), (y_1, y_2) \ra =
\ll x_1, y_1 \gg -
\ll x_2, y_2 \gg, \quad x_1, x_2, y_1, y_2 \in \g,
\]
where $\ll.,.\gg$ is a nondegenerate invariant symmetric bilinear 
form on $\g$ whose restriction to a compact real form of $\g$ is 
negative definite. Lagrangian splittings of $(\gog, \lara)$ have
been classified by Delorme \cite{delorme:manin-triples}. In particular, one 
has the Belavin--Drinfeld splittings $\gog = \g_{{\rm diag}} + \l$, where
$\g_{{\rm diag}}$ is the diagonal of $\gog$. 
For any $\l \in \L(\gog)$, $N(\l)$-orbits in $\L(\gog)$ 
can be described by using
results in \cite{LuY}. For an 
arbitrary Lagrangian splitting $\gog = \l_1 + \l_2$
we prove that the intersection of each $N(\l_1)$ and
$N(\l_2)$-orbit on $\L(\gog)$ is connected. Further,
using \cite{LuY}, we compute the rank of 
all corresponding Poisson structures 
$\Pi_{\l_1, \l_2}$ on $\L(\gog)$. This result extends
the dimension formulas for symplectic leaves in 
the second author's classification \cite{Y} of symplectic 
leaves of Belavin--Drinfeld Poisson structures on 
complex reductive Lie groups. Our result also
generalizes the rank formulas of S. Evens and the first
author for the standard Poisson structure on $\Lgg$.   

As have been shown in \cite{E-Lu1, E-Lu2}, all real and complex semi-simple
symmetric spaces, as well as certain of their compactifications
can be embedded into suitable varieties of Lagrangian subalgebras. Out 
results
show that all such spaces carry Poisson structures and natural partitions
into regular Poisson subvarieties. 

All manifolds and vector spaces in this paper, unless otherwise stated, 
are assumed to be either complex or real. 

A submanifold $N$ of a 
Poisson manifold $(M, \pi)$ will be called a {\em{complete Poisson
submanifold}} if it is closed under all Hamiltonian
flows or equivalently it is a union of symplectic leaves
of $\pi$.

{\bf{Acknowledgements.}} The second author would like to thank the 
University
of Hong Kong for the warm hospitality during his visits in March 2004 and 
August 2005 when this work was initiated. The first author 
would like to thank UC Santa Barbara for her visit in July 2006 during
which the paper took its final form.
We would also like to thank Xuhua He for a key argument in the proof of 
\prref{NN-connected}.
The first author was partially 
supported by HKRGC grants 703304 and 703405, and the second author by
NSF grant DMS-0406057 and an Alfred P. Sloan research fellowship.

\sectionnew{The Poisson spaces $D/Q$}\lb{sec-DQ}

Recall that a quadratic Lie algebra is a pair $(\d, \lara)$, where
$\d$ is a Lie algebra and $\lara$ is an invariant symmetric 
nondegenerate bilinear form on $\d$.
Throughout this section, we fix a quadratic Lie algebra $(\d, \lara)$
and a connected Lie group $D$ with Lie algebra $\d$. For a subspace
$V$ of $\d$, set 
\begin{equation}
\label{ort}
V^\perp = \{x \in \d \mid \la x, \, y \ra = 0, \forall
y \in V \}.
\end{equation}

\subsection{Lagrangian splittings}

A {\it coisotropic} (resp. 
{\it Lagrangian, isotropic}) subalgebra of $\d$ (with respect to $\lara)$
is by definition 
a Lie subalgebra 
$\q$ of $\d$ such that  $\q^\perp \subset \q$ (resp. 
$\q^\perp = \q, \; \q \subset \q^\perp$).

\bde{de-lagr-splitting}
A Lagrangian splitting of $\d$ is a 
vector space direct sum decomposition $\d = \u + \up$, where
$\u$ and $\up$ are both Lagrangian subalgebras of $\d$.
The triple $(\d, \u, \up)$ 
is also called a Manin triple \cite{k-s:quantum}.
\ede
Given a Lagrangian splitting $\d = \u + \up$, for a subspace $W \subset \u$,
set
\begin{equation}
\label{0comp}
W^0 = \{\xi \in \up \mid \la \xi, \, x \ra = 0, \forall
x \in W \} = W^\perp \cap \up.
\end{equation}

We now recall how Lagrangian splittings give rise to Poisson Lie groups.
Recall that a Poisson Lie group is a pair $(G, \pi)$, where $G$
is a Lie group and $\pi$ is a Poisson structure
on $G$ such that the group multiplication $G \times G \to G$ is a 
Poisson map. When a (not necessarily closed) subgroup $H$ of $G$ is also a 
Poisson submanifold with
respect to $\pi$, $(H, \pi)$ is itself a Poisson Lie group and is called a 
Poisson Lie subgroup
of $(G, \pi)$.
If $(G, \pi)$ is a Poisson Lie group, then $\pi(e) = 0$, where 
$e \in G$ is the identity element. Let $\g$ be the 
Lie algebra of $G$, and let $d_e \pi: \g \to \wedge^2 \g$
be the linearization of $\pi$ at $e$ defined by
\[
(d_e \pi)(x) = (L_{\widetilde{x}} \pi)(e),
\]
where for $x \in \g, \widetilde{x}$ is any 
local vector fields with $\widetilde{x}(e) = x$
and $L_{\tilde{x}}\pi$ is the Lie derivative of $\pi$ at $e$.
Then $(\g, d_e \pi)$ is a Lie bialgebra
\cite{k-s:quantum} called the tangential Lie bialgebra of $(G, \pi)$.

Assume that $\d = \u + \up$ is a Lagrangian splitting. The bilinear
form $\lara$ induces a non-degenerate pairing between $\u$ and $\up$. 
Define
\begin{align}
\lb{deltau}
& \delta_{\u}: \; \u \lra \wedge^2 \u: \; 
\la \delta_\u (x), \, y \wedge z \ra = 
\la x, [y, z]\ra, \hs & x \in \u, y, z \in \up,\\
\lb{deltaup}
& \delta_{\up}: \; \up \lra \wedge^2 \up: \; 
\la \delta_{\up}(x), y \wedge z \ra
=\la x, [y, z]\ra, \hs & x \in \up, y, z \in \u.
\end{align}
Then $(\u, \delta_\u)$ and $(\up, \delta_{\up})$ are Lie bialgebras 
\cite{k-s:quantum}.
Associated to the splitting $\d = \u + \up$ we also have the $r$-matrix
\begin{equation}\lb{ruu}
R_{\u,\up} = \frac{1}{2} \sum_{j=1}^{n} \xi_j \wedge x_j \in \wedge^2 \d,
\end{equation}
where $\{x_1, x_2, \cdots, x_n\}$ and 
$\{\xi_1, \xi_2, \cdots, \xi_n\}$ are  bases
of $\u$ and $\up$, respectively, such that 
$\la x_i, \xi_j \ra = \delta_{ij}$
for $1 \leq i, j \leq n$.
It is easy to see that 
$R_{\u, \up}$ is independent of the 
choice of the bases. Moreover, the Schouten bracket 
$[R_{\u, \up}, \,R_{\u, \up}] \in \wedge^3 \d$ is given by
\begin{equation} \lb{Sch_r}
\la 
[R_{\u, \up}, \,R_{\u, \up}], \, a \wedge b \wedge c \ra = 
2\la a, [b, c] \ra, \hs a, b, c \in \d.
\end{equation}
Recall that $D$ is a connected Lie group with Lie algebra $\d$.
Denote by $U$ and $\Up$
the connected subgroups of $D$ with Lie algebras $\u$ and $\up$, 
respectively. Let
$R_{\u, \up}^{l}$ and $R_{\u, \up}^r$ be the left and the right 
invariant bi-vector fields on $D$ with values $R_{\u, \up}$ at
the identity element.  Set
\begin{equation}\lb{piDuup}
\pi_{\u, \up}^{D} := R_{\u, \up}^r - R_{\u, \up}^l.
\end{equation}
The following fact can be found in \cite{ES, k-s:quantum}.

\bpr{double}
The bivector field $\pi_{\u, \up}^{D}$
is a Poisson structure on $D$ and $(D, \pi^D_{\u, \up})$ is a Poisson Lie 
group.
Both $U$ and $\Up$ are Poisson Lie subgroups of $(D, \pi^D_{\u, \up})$. 
Let 
\begin{equation}\lb{piU-piUp}
\pi_U = \pi^D_{\u, \up}|_U, \hs \pi_{\Up} = -\pi^D_{\u, \up}|_{\Up}.
\end{equation}  
Then the tangential Lie bialgebras of the Poisson Lie groups
$(U, \pi_U)$ and $(\Up, \pi_{\Up})$ 
are respectively $(\u, \delta_\u)$ and $(\up, \delta_{\up})$.
\epr

\subsection{The Poisson spaces $D/Q$}
\label{sec_DQ}
Assume that $Q$ is 
a closed subgroup of $D$ whose Lie algebra $\q$ is a coisotropic
subalgebra of $(\d, \lara)$.
For an integer $k \geq 1$,
let $\chi^k(D/Q)$ be the space of $k$-vector fields on $D/Q$. 
Then the left action of $D$ on $D/Q$ gives rise to 
the Lie algebra anti-homomorphism 
\[
\kappa \colon \; \d \lra \chi^1(D/Q)
\]
whose multi-linear extension $\wedge^k \d \to \chi^k(D/Q)$ will
be denoted by the same letter. 

Given a Lagrangian splitting $\d = \u + \up$, define
the bivector field  $\Pi_{\u, \up}$ on $D/Q$ by
\begin{equation}
\lb{piuu}
\Pi_{\u, \up} :=\kappa(R_{\u, \up}),
\end{equation}
recall \eqref{ruu}. The following theorem
is the main result for this subsection.

\bth{kappaR} For every 
Lagrangian splitting $\d = \u + \up$ and 
every closed subgroup $Q$ of $D$ whose Lie algebra $\q$ is
coisotropic in $\d$,

1) $\Pi_{\u, \up}$ is a Poisson bi-vector field on $D/Q$;

2) all $U$ and $U^\prime$-orbits in $D/Q$ are 
complete Poisson submanifolds of $(D/Q, \Pi_{\u, \up})$.
\eth

\begin{proof} 1) The Lie algebra of the stabilizer subgroup of each 
point of $D/Q$ is a coisotropic subalgebra of $\d$. To prove 
that $\Pi_{\u, \up}$ is Poisson, 
it suffices to show that 
\[
[R_{\u, \up}, \, R_{\u, \up}] \in \q \wedge \d \wedge \d
\] 
for each coisotropic subalgebra $\q$ of $\d$. This is 
equivalent to 
\[
\la [R_{\u, \up},\, R_{\u, \up}], a \wedge b \wedge c \ra = 0, 
\quad \forall a, b, c \in \q^\perp
\] 
which follows from \eqref{Sch_r} because 
$\q^\perp$ is an isotropic subalgebra of $\d.$ 

Denote by $\kappa_\q: \d \to \d /\q$ the
canonical projection and its induced map
$\wedge^2 \d \to \wedge^2 (\d / \q)$. The second part 
of \thref{kappaR} now
follows from \leref{lem_q-1} below.
\end{proof}

\ble{lem_q-1} 
For every coisotropic subalgebra $\q$ of $\d$, one has
\[
\kappa_\q   (R_{\u, \u'}) \in \left(\kappa_\q (\wedge^2 \u)\right) 
\cap \left(\kappa_\q(\wedge^2 \up)\right).
\]
\ele

\begin{proof}
It is sufficient to show that 
$\kappa_\q (R_{\u, \up}) \in 
\kappa_\q (\wedge^2 \u)$. Let $\{x_1, x_2, \cdots, x_l\}$ 
be a basis for $\u \cap \q$. Extend it to a basis 
$\{x_1, x_2, \cdots, x_l, x_{l+1}, \cdots, x_n\}$ of $\u$.
Let $\{\xi_1, \cdots, \xi_n\}$ be the dual basis of $\up$
with respect to $\lara$.
It is easy to see that $(\u \cap \q)^0 = 
{\rm Span}\{\xi_{l+1}, \xi_{l+1}, \cdots, \xi_n\} =
p_{\up}(\q^\perp)$, 
recall \eqref{0comp}, where
$p_{\up}: \d \to \up$ is the projection along $\u$.
Choose $y_j \in \u$
such that 
$y_j+ \xi_j \in \q^\perp$ for $l+1 \leq j \leq n$ 
and write
\[
R_{\u, \u'} = \frac{1}{2} 
\sum_{j=1}^{l} \xi_j \wedge x_j + \frac{1}{2} 
\sum_{j=l+1}^{n} (y_j +\xi_j) \wedge x_j -
\frac{1}{2} \sum_{j=l+1}^{n} y_j \wedge x_j.
\]
Since $\kappa_\q(x_j) = 0$ for $1 \leq j \leq l$ and
$\kappa_\q( y_j+\xi_j) =0$
for $l+1 \leq j \leq n$, we have
\begin{equation}\lb{kappa-Ruu}
\kappa_\q (R_{\u, \u'}) =  
- \frac{1}{2} \sum_{j=l+1}^{n} \kappa_\q(y_j) 
\wedge \kappa_\q(x_j) \in \kappa_\q  (\wedge^2 \u).
\end{equation}
\end{proof}
 
In the setting of 
\thref{kappaR}, the action map 
\[
(D, \pi^D_{\u, \up})  \times (D/Q, \, \Pi_{\u, \up}) 
\lra (D/Q, \, \Pi_{\u, \up})
\]
is easily seen to be Poisson. Thus $(D/Q, \Pi_{\u, \up})$ is a Poisson 
homogeneous space \cite{Drinfeld} of 
$(D, \pi^D_{\u , \up})$.
By part 2) of \thref{kappaR}, 
each $U$ and $\Up$-orbit in $D/Q$
is a Poisson homogeneous space of $(U, \pi_U)$ and
$(\Up, -\pi_{\Up})$, respectively.
\subsection{Rank of the Poisson structure $\Pi_{\u, \u^\prime}$ on $D/Q$}
For $d \in D$ set $\ud = dQ \in D/Q$.
Consider the $U$-orbit $U . \ud$ 
through $\ud$. Then $(U . \ud, \Pi_{\u, \up})$ is a Poisson homogeneous 
space of $(U, \pi_U)$. Denote by $\l_{\ud}$ the {\it Drinfeld Lagrangian subalgebra} of $\d$,
associated to the base point $\ud$ of $U_\cdot \ud$, 
cf. \cite{Drinfeld}. It is defined as follows: 
identify $T_{\ud} (U_\cdot \ud) \cong \u /(\u \cap \Ad_d q)$
and regard $\Pi_{\u, \up}(\ud)$ as an element in 
$\wedge^2 \left(\u/(\u \cap \Ad_d q)\right)$.
For $\xi \in (\u \cap \Ad_d q)^0$, recall \eqref{0comp}, let
$\iota_\xi \Pi_{\u, \up}(\ud) \in \u /(\u \cap \Ad_d q)$ be such that
$\la \iota_\xi \Pi_{\u, \up}(\ud), \; 
\eta \ra = \Pi_{\u, \up}(\ud)(\xi, \eta)$
for all $\eta \in (\u \cap \Ad_d q)^0$. Then $\l_{\ud} \subset \d$ is given by
\[
\l_{\ud} = \{x + \xi \mid x \in \u, \xi \in (\u \cap \Ad_d q)^0, 
\,\iota_\xi \Pi_{\u, \up}(\ud) = x + 
\u \cap \Ad_d q\}.
\]
If $\Rank_{\Pi_{\u, \u^\prime}}(\ud)$ denotes 
the rank of $\Pi_{\u, \u^\prime}$ at $\ud$, it is easy to see from the
definition of $\l_{\ud}$ that 
\begin{equation}\lb{rank-ld}
\Rank_{\Pi_{\u, \u^\prime}}(\ud) = 
\dim (U_\cdot \ud) - \dim( \up \cap \l_{\ud}).
\end{equation}

\bpr{ld} For any $d \in D$, the Drinfeld Lagrangian subalgebra $\l_{\ud}$ is
\begin{equation}\lb{eq-ld}
\l_{\ud} = \Ad_d \q^\perp + \u \cap \Ad_d q.
\end{equation}
\epr

\begin{proof}
Since the stabilizer subalgebra of $\d$ at $\ud$ is $\Ad_d \q$, it is
enough to prove that $\l_{\underline{e}} = \q^\perp + \u \cap \q$, 
where $e$ is the identity element of $D$.
Note that since $\q^\perp \subset \q$,
\[
(\q^\perp + \u \cap \q)^\perp = 
\q \cap (\u \cap \q)^\perp = \q \cap (\u + \q^\perp) = 
\u \cap \q + \q^\perp.
\]
Thus 
$\q^\perp + \u \cap \q$ is a Lagrangian subspace of $\d$. 
Since $\l_{\underline{e}}$ is Lagrangian in $\d$ and 
$\u \cap \q \subset \l_{\underline{e}}$, it is sufficient to show that
$\q^\perp \subset \l_{\underline{e}}$.
 
For $1 \leq i \leq n$ and $l+1 \leq j \leq n$, 
let $x_i \in \u, \xi_i \in \up$ and $y_j \in \u$
be as in the proof of \leref{lem_q-1}. 
If $y + \xi \in \q^\perp$, for some
$y \in \u, \xi \in \up$, then
$\xi = \sum_{j=l+1}^{n} \lambda_j \xi_j$. Thus 
$y + \xi - \sum_{j=l+1}^{n} \lambda_j (y_j+\xi_j) 
\in \u \cap \q^\perp \subset \l$.
The proposition will now follow if we show 
that $y_j + \xi_j \in \l_{\underline{e}}$ for every $l+1 \leq j \leq n$.
By \eqref{kappa-Ruu}
\[
\Pi_{\u, \u^\prime}(\ue) = -\frac{1}{2} 
\sum_{j=l+1}^{n} (y_j + \u \cap \q) \wedge (x_j+\u \cap \q)
\in \wedge^2 (\u /(\u \cap \q)) \cong \wedge^2 T_{\ue} D/Q.
\]
Thus for each $l+1 \leq j \leq n$,
\[
\iota_{\xi_j} \Pi_{\u, \u^\prime}(\ue) = \frac{1}{2} y_j -\frac{1}{2} 
\sum_{k=l+1}^{n} \la \xi_j, y_k \ra x_k + \u \cap \q.
\]
Since $0 = \la y_j + \xi_j, y_k + \xi_k \ra = 
\la \xi_j, y_k \ra + \la \xi_k, y_j \ra$ for 
$l+1 \leq k \leq n$ and since $x_j \in \u \cap \q$ for $1 \leq j \leq l$, we 
have 
\begin{align*}
\iota_{\xi_j} \Pi_{\u, \u^\prime}(\ue) &= 
\frac{1}{2} y_j +\frac{1}{2} \sum_{k=l+1}^{n} 
\la \xi_k, y_j \ra x_k + \u \cap \q 
= \frac{1}{2} y_j +\frac{1}{2} \sum_{k=1}^{n} 
\la \xi_k, y_j \ra x_k + \u \cap \q \\
&= \frac{1}{2} y_j + \frac{1}{2} y_j + \u \cap \q
=   y_j + \u \cap \q,
\end{align*}
we see that $y_j + \xi_j \in \l_{\underline{e}}$. 
\end{proof}

Since $\u + \up = \d$,  $U$-orbits and $U^\prime$-orbits in $D/Q$ intersect
transversally. By \thref{kappaR}, any such non-empty intersection is a Poisson 
submanifold of $\Pi_{\u, \up}$. The following corollary gives the corank of
$\Pi_{\u, \up}$ in $U_\cdot \ud \cap U^\prime_\cdot \ud$ 
at $\ud$ for every $d \in D$.

\bco{corank-UU} For any $d \in D$,
\[
\Rank_{\Pi_{\u, \up}}(\ud) = \dim (U_\cdot \ud \cap U^\prime_\cdot \ud)+
\dim (D/Q)
-\dim (U^\prime_\cdot \ud) -\dim(\up \cap \l_{\ud}),
\]
where $\l_{\ud}$ is the Drinfeld Lagrangian subalgebra given by \eqref{eq-ld}.
\eco

\begin{proof} The statement follows immediately from \eqref{rank-ld}
and the fact that 
\[
\dim(U_\cdot \ud) + \dim(U^\prime_\cdot \ud) - \dim(D/Q) = 
\dim(U_\cdot \ud \cap U^\prime_\cdot \ud).
\]
\end{proof}

\subsection{First main theorem}\lb{first-main}
Recall that a manifold with a Poisson structure of constant rank
 is called a 
regular Poisson manifold.

\bth{qqqperp} 
If $\q$ is a coisotropic subalgebra 
of $\d$ such that $[\q, \q] \subset \q^\perp$, 
then for any closed subgroup $Q$ of $D$ 
with Lie algebra $\q$ and 
for any Lagrangian splitting $\d = \u + \up$ of $\d$, the intersection of
any $U$-orbit with any $U^\prime$-orbit in $D/Q$ is 
a regular Poisson submanifold for 
the Poisson structure $\Pi_{\u, \up}$.
\eth

\begin{proof}
Let again $e$ be the identity in $D$.
Since $[\Ad_d\q, \Ad_d\q] \subset (\Ad_d\q)^\perp$ 
for any $d \in D$, it suffices to show 
that $U_\cdot \ue \cap U^\prime_\cdot \ue$ is a 
regular Poisson manifold of $\Pi_{\u, \up}$.
Let $a \in U$ and $a^\prime \in U^\prime$ be such 
that $\ua = \ua^\prime 
\in U_\cdot \ue \cap U^\prime_\cdot \ue$. Then there 
exists $b \in Q$ such that $a = a^\prime b$. 
Thus
\begin{align*}
\dim(\up \cap \l_{\underline{a}}) &=\dim(\up \cap \Ad_a(\q^\perp + \u \cap \q)) = 
\dim(\up \cap \Ad_{a^\prime b}(\q^\perp + \u \cap \q)) \\
&=\dim(\up \cap \Ad_b(\q^\perp + \u \cap \q)).
\end{align*}
Since $[\q, \q] \subset \q^\perp$, we have $[\q, \q^\perp + \u \cap \q] 
\subset [\q, \q] \subset \q^\perp \subset \q^\perp + \u \cap \q$. Thus 
$Q$ normalizes $\q^\perp + \u \cap \q$, and
$\Ad_b(\q^\perp + \u \cap \q) = \q^\perp + \u \cap \q$. Hence
$\dim(\up \cap \l_{\underline{a}}) = \dim(\up \cap \l_{\underline{e}})$. It 
follows from \coref{corank-UU} that
the rank of $\Pi_{\u, \u^\prime}$ at $\underline{a}$ 
is the same as that at $\ue$.
\end{proof}

The special case of \thref{qqqperp}, 
when the Drinfeld subalgebras of all points of 
$(D/Q, \Pi_{\u, \up})$ integrate to closed subgroups of $D$, 
can be also proved using Karolinsky's result \cite{karo:sl}.

\bre{re-qqq}
Note that if $\q \subset \d$ is a coisotropic subalgebra 
such that $\n(\q) = \q$, where
$\n(\q)$ is the normalizer of $\q$ in $\d$, 
then $[\q, \q] \supset \q^\perp$.
Indeed, if $x \in [\q, \q]^\perp$, then
$\la x, [\q, \q] \ra = 0$ which implies that $\la [x, \q], \q \ra = 0$, so 
$[x, \q] \subset \q^\perp \subset \q$. Thus $x \in \n(\q) = \q$. This shows 
that $[\q, \q]^\perp \subset \q$, so $[\q, \q] \supset \q^\perp$. We thus 
conclude that, if $\q$ is coisotropic such that $\n(\q)=\q$ and 
$[\q, \q] \subset \q^\perp$, then $[\q, \q] = \q^\perp$. 
\ere

This remark will be used in $\S$\ref{second-main}.

\bco{DL}
Let $L$ be a closed subgroup of $D$ whose Lie algebra $\l \subset \d$ is
Lagrangian.  Then for any Lagrangian
splitting $\d = \u + \up$,  symplectic leaves of $\Pi_{\u,\up}$ in $D/L$
are precisely the connected components of the intersections of
$U$ and $U^\prime$-orbits in $D/L$.
\eco

\begin{proof}
Again it is enough to prove that  
$U_\cdot \ue \cap U^\prime_\cdot \ue$ is symplectic.
By \coref{corank-UU} and \thref{qqqperp}, the corank of $\Pi_{\u, \up}$ in 
$U_\cdot \ue \cap U^\prime_\cdot \ue$ is equal to 
\[
\dim(U^\prime_\cdot \ue) + \dim(\up \cap \l) - \dim(D/L) = 0.
\]
\end{proof}

\bre{homog} 
\coref{DL} describes the symplectic leaves for a large class of 
Poisson homogeneous spaces. Indeed, by \cite{Drinfeld}, every
Poisson homogeneous space of $(U, \pi_U)$ is of the form 
$U/H$, where $H$ is a subgroup of $U$ whose Lie algebra is 
$\u \cap \l$ for a Lagrangian subalgebra $\l$ of $\d$. In the
case when $H = U \cap L$, where $L$ is a closed subgroup of $D$ with Lie 
algebra $\l$, we can embed $U/H$ into $D/L$ as the $U$-orbit through $\ue$. 
This embedding is also Poisson. Thus the symplectic leaves of $U/H$ are the 
connected components of the intersections of $U/H$ with the $U^\prime$-orbits 
in $D/L$, i.e. with the $(U^\prime, -\pi_{U^\prime})$ Poisson homogeneous 
spaces inside $(D/L, \Pi_{\u, \u^\prime})$.
\ere

We conclude this subsection with an example showing that 
the statement of \thref{qqqperp} is incorrect if the 
condition $[\q, \q] \subset \q^\perp$ is dropped.

\bex{counterex} Let $G$ a connected complex simple Lie group
with a pair of opposite Borel subgroups $B$ and $B^-$.
Set $\g = \Lie G$, $T= B \cap B^-$, and $\h = \Lie T$. 
Then $\d = \g \oplus \h$ is a quadratic Lie algebra with the 
bilinear form
\[
\la (x_1,x_2),(y_1,y_2)\ra =\ll x_1,y_1\gg- \ll x_2,y_2\gg, 
\hspace{.2in} x_1, y_1, 
\in \g, x_2, y_2 \in \h
\]
where $\ll.,.\gg$ is a nondegenerate symmetric invariant bilinear form on 
$\g$. Let $D = G \times T$. 
Given a parabolic subgroup $P \supset B$ of $G$, the Lie algebra
$\q$ of $Q=P \times T$ is coisotropic but does not satisfy
the condition $\q^\perp \subset [\q, \q]$. The following subalgebras
provide a Lagrangian splitting of $\g \oplus \h$:
\begin{equation*}
\u = \{(x+h, h) \mid x \in \n, \, h \in \h \}, \hspace{.2in}
\up = \{(x+h, -h) \mid x \in \n^-, \, h \in \h \}, 
\end{equation*}
where $\n$ and $\n^-$ are 
the nilpotent radicals of $\Lie B$ and $\Lie B^-$.
Under the identification $(G \times T)/ (P \times T) \cong G/P$ the 
Poisson structure $\Pi_{\u, \up}$ corresponds to the Poisson 
structure 
\[
\pi= \kappa \left( \sum_{\alpha \in \Delta^+} 
f_\alpha \wedge e_\alpha \right),
\]
where, $\Delta^+$ is the set of positive roots of $\g$ corresponding
to $\n$,
$\{ e_\alpha \}$ and $\{ f_\alpha \}$ are sets of  root vectors of 
$\g$, normalized by $\ll e_\alpha, f_\alpha \gg = 1$,
and $\kappa$ is the extension to $\wedge^2 \g$ of the infinitesimal action 
of $\g$ on $G/P$. It was shown in \cite{GY} that the partition of 
$(G/P, \pi)$ by $T$-orbits of leaves (which is a partition by regular 
Poisson submanifolds) coincides with Lusztig's partition \cite{Lus} of 
$G/P$. The strata of this partition are 
\begin{equation}
\label{Lusztig-partition}
{\mathrm{pr}}_P \left( B w_1 . B \cap B^- w_2 . B \right), \quad
w_1 \in W, w_2 \in W^{W_P}_{\mathrm{max}}.  
\end{equation}
Here ${\mathrm{pr}}_P : \; G/B \to G/P$ denotes the standard projection, 
$W$ the Weyl group of $(G,T)$, and $W^{W_P}_{\mathrm{max}}$ 
the set of the maximal length representatives of cosets in $W/W_P$ where,
$W_P$ is the parabolic subgroup of $W$ corresponding to $P$.

Under the identification $(G \times T) / (P \times T) \cong G/P$ the $U$ and 
$\Up$-orbits on $(G \times T) /( P \times T)$ correspond respectively to the 
$B$ and $B^-$-orbits on $G/P$. The coarser partition of 
\eqref{Lusztig-partition} 
by intersecting $B$ and $B^-$-orbits on $(G/P, \pi)$ is no longer a partition
by regular Poisson submanifolds if $P \neq B$ or $G$. This is easily seen by 
applying \thref{rank-main} below and the Poisson embedding 
\cite[(1.10)]{GY}.
\eex
\subsection{Intersections of $N(\u)$ and $N(\up)$-orbits}\lb{NuNup-orbits}
For a Lagrangian splitting $\d = \u + \up$, let $N(\u)$ and $N(\up)$ be the
normalizer subgroups of $\u$ and $\up$ in $D$, respectively. Both
$N(\u)$ and $N(\up)$ are closed subgroups of $D$, and sometimes $N(\u)$
and $N(\up)$-orbits in a space $D/Q$ are easier to determine than
the $U$ and $\Up$-orbits. This is the case for the examples considered in
$\S$\ref{Lgg}. In this subsection, we prove some facts on $N(\u)$ and
$N(\up)$-orbits.

It is clear from \thref{kappaR} that for any closed subgroup $Q$ of $D$ with
coisotropic Lie subalgebra $\q$ and for any Lagrangian splitting 
$\d = \u + \up$, 
all $N(\u)$ and $N(\up)$-orbits in $D/Q$ are complete Poisson 
submanifolds with respect to the Poisson structure $\Pi_{\u, \up}$,
cf. \thref{kappaR}.
Recall the Poisson structure $\pi^D_{\u, \up}$ on $D$ from \eqref{piDuup}. 

\ble{NuNup}
For any Lagrangian splitting $\d = \u + \up$, 
the Poisson structure 
$\pi^D_{\u, \up}$ on $D$ vanishes at all points in $N(\u) \cap N(\up)$. 
Consequently, for any closed subgroup $Q$ of $D$ with coisotropic Lie 
subalgebra $\q$, 
$N(\u) \cap N(\up)$ leaves the Poisson structure 
$\Pi_{\u, \up}$ on $D/Q$ invariant.
\ele

\begin{proof} Let $d \in N(\u) \cap N(\up)$. If $\{x_1, \ldots, x_n\}$ and 
$\{\xi_1, \ldots, \xi_n\}$ is a pair of dual bases for $\u$ and $\up$ with 
respect to $\lara$, then so are $\{\Ad_d(x_1), \ldots \Ad_d(x_n)\}$ and 
$\{\Ad_d(\xi_1), \ldots, \Ad_d(x_n)\}$. Thus
\[
\pi^D_{\u, \up}(d) = L_d(r_{\u, \up}) - R_d(r_{\u, \up})=
R_d\left(
\sum_{i=1}^n \Ad_d(\xi_i) \wedge \Ad_d(x_i)  
-\sum_{i=1}^n \xi_i \wedge x_i
\right)=0.
\]
Since
the $(D, \pi^D_{\u, \up})$-action on $(D/Q, \Pi_{\u, \up})$ is
Poisson,
$N(\u) \cap N(\up)$ leaves $\Pi_{\u, \up}$ invariant.
\end{proof}

Next we generalize \thref{qqqperp} to intersections
of arbitrary $N(\u)$-orbits and $N(\up)$-orbits in $D/Q$.
Its proof is similar to the one of \thref{qqqperp}
and is left to the reader.

\bpr{N-regular}
If $\q$ is a coisotropic subalgebra 
of $\d$ such that $[\q, \q] \subset \q^\perp$, then for any closed subgroup 
$Q$ with Lie algebra $\q$ and 
for any Lagrangian splitting $\d = \u + \up$ of $\d$, the intersection of
any $N(\u)$-orbit with any $N(\up)$-orbit in $D/Q$ 
is a regular Poisson submanifold for 
the Poisson structure $\Pi_{\u, \up}$.
\epr

Although \prref{N-regular} provides a stronger 
result than \thref{qqqperp}, it is apriori possible that 
the geometry of the strata of the coarser partition from 
\prref{N-regular}
is more complicated than that of the strata of the finer
partition from \thref{qqqperp}. The next result, \prref{orbit},
shows that this
is not the case. First we prove an auxiliary lemma.

\ble{Nu}
Let $\d = \u + \up$ be any Lagrangian splitting of $\d$. Assume that 
$N(\u)$ is connected.  Then $N(\u) = (U^\prime \cap N(\u))^o U$, where
$(U^\prime \cap N(\u))^o$ denotes the identity component of
the group $U^\prime \cap N(\u)$. Moreover, $N(\u)$ is a Poisson
Lie subgroup of $(D, \pi)$.
\ele
 
\begin{proof}
This is because $(U^\prime \cap N(\u))^o U$ is a connected subgroup of
$D$ with Lie algebra $\up \cap \n(\u) + \u$ which is 
equal to $\n(\u)$ because $\d=\u + \up$.
\end{proof}

\bpr{orbit}
Let $\d = \u + \up$ be any Lagrangian splitting of $\d$, and assume that 
$N(\u)$ and $N(\up)$ are both connected. Let $X$ be any Poisson space with 
a Poisson $(D, \pi^D_{\u, \up})$-action. 
Let $x \in X$ such that $N(\u)x \cap N(\up)x \neq \emptyset$.   
Then the group $N(\u) \cap N(\up)$ acts transitively on the
set of intersections of $U$-orbits and $U^\prime$-orbits in 
$N(\u)x \cap N(\up)x$.
\epr
\begin{proof}
Using \leref{Nu} we obtain
\begin{align*}
N(\u)x \cap N(\up)x & = \bigcup_{\alpha \in (U^\prime \cap N(\u))^o, 
\beta \in (U  \cap N(\up))^o} (\alpha Ux) \cap ( U^\prime \beta x)\\
& = \bigcup_{\alpha \in (U^\prime \cap N(\u))^o, 
\beta \in (U  \cap N(\up))^o} \alpha(Ux \cap U^\prime \beta x) \\&= 
\bigcup_{\alpha \in (U^\prime \cap N(\u))^o, 
\beta \in (U  \cap N(\up))^o} \alpha \beta (Ux \cap U^\prime x).
\end{align*}

\end{proof}

We finish this section with a formula for 
the corank of $\Pi_{\u, \up}$ to be used in
$\S$\ref{Lgg}. As before $Q \subset D$ 
is assumed to be a closed subgroup with coisotropic Lie 
subalgebra $\q$.
For any $d \in D$, let 
\[
{\rm Corank}_{\Pi_{\u, \up}}(\ud) = 
\dim (N(\u)_\cdot \ud \cap N(\up)^\prime_\cdot \ud)-
\Rank_{\Pi_{\u, \up}}(\ud)
\]
be the corank of $\Pi_{\u, \up}$ in $N(\u)_\cdot \ud \cap N(\up)_\cdot \ud$ 
at $\ud \in D/Q$.

\ble{formula-NN} In the above setting,
for any $d \in D$,
\begin{align*}
{\rm Corank}_{\Pi_{\u, \up}}(\ud) &= 
\dim \n(\up)-\dim(D/Q) +\dim \n(\u) -\dim \u \\
&-\dim (\n(\u) \cap \Ad_d \q) + \dim (\u \cap \Ad_d \q)\\
&-\dim(\n(\up) \cap \Ad_d\q) + \dim(\up \cap \l_{\ud}),
\end{align*}
where $\l_{\ud}=\Ad_d \q^\perp + \u \cap \Ad_d q$ is the 
Drinfeld Lagrangian subalgebra as in 
\eqref{eq-ld}.
\ele

\begin{proof} Since $N(\u)_\cdot \ud$ and $N(\up)_\cdot \ud$ 
intersect transversally, 
\begin{align*}
\dim (N(\u)_\cdot \ud \cap N(\up)^\prime_\cdot \ud) &= 
\dim(N(\u)_\cdot \ud) +
\dim(N(\up)_\cdot \ud) - \dim(D/Q) \\
&= \dim \n(\u) - \dim(\n(\u) \cap \Ad_d \q) \\
&+ \dim \n(\up) - \dim(\n(\up) \cap \Ad_d \q) - \dim(D/Q).
\end{align*}
By \eqref{rank-ld}, $\Rank_{\Pi_{\u, \up}}(\d) = 
\dim \u - \dim (\u \cap \Ad_d \q) - \dim(\up \cap \l_{\ud}).$
The formula for ${\rm Corank}_{\Pi_{\u, \up}}(\ud)$ in 
\leref{formula-NN} thus follows.
\end{proof}
\sectionnew{The variety of Lagrangian subalgebras associated to 
a reductive Lie algebra}\lb{Lag}
\subsection{General case}\lb{general}
Let $(\d, \lara)$ be a $2n$-dimensional quadratic Lie algebra and let $D$
be a connected Lie group with Lie algebra $\d$. We will denote by $\Ld$
the variety of all Lagrangian subalgebras of $\d$. It is an algebraic 
subvariety of the Grassmannian $\Gr(n, \d)$ of $n$-dimensional 
subspaces of $\d$. The group $D$ acts on $\Ld$ through the
adjoint action. 

Fix a Lagrangian splitting $\d = \u + \up$, recall
$R_{\u, \up} \in \wedge^2 \d$ given by \eqref{ruu}.
Let again
$\kappa: \d \to \chi^1(\Ld)$ be the Lie algebra 
anti-homomorphism from $\d$ to the
Lie algebra of vector fields on $\Ld$, and define the bi-vector field
$\Pi_{\u, \up} = \kappa(R_{\u, \up})$ on $\Ld$.
For $l \in \Ld$, let $N(\l)$ and $\n(\l)$ be 
respectively the normalizer subgroup of $\l$
in $D$ 
and the normalizer subalgebra of $\l$ in $\d$.
Then the $D$-orbit in $\Ld$ through $\l$ is isomorphic to $D/N(\l)$.
Clearly,
$\n(\l)$ is coisotropic in $\d$ because it contains $\l$. Thus it
follows from  \thref{kappaR} that $\Pi_{\u, \up}$ is a 
Poisson structure on $\Ld$, see also \cite{E-Lu1}. 
The following proposition now follows immediately from  
\thref{qqqperp}.

\bpr{Ld-regular-1} Assume that $(\d, \lara)$ is an even dimensional 
quadratic Lie algebra and $D$ is a connected
Lie group with Lie algebra $\d$ such that 
for every  $\l \in \Ld$
\[
[\n(\l), \n(\l)] \subset (\n(\l))^\perp.
\]
Then for any Lagrangian splitting
$\d = \u + \up$, the intersection of an $N(\u)$-orbit and 
an $N(\up)$-orbit in $\Ld$ is a regular Poisson submanifold for
the Poisson structure $\Pi_{\u, \up}$.
\epr

\subsection{Second main theorem: the case of a 
reductive Lie algebra}\lb{second-main}

When $(\d, \lara)$ is a reductive quadratic Lie algebra, we have the 
following second main theorem of the paper.

\bth{Ld-regular-2} If $D$ is a connected
complex or real reductive Lie group and
$\lara$ is a nondegenerate symmetric invariant bilinear form on $\d = \Lie D$, 
then for any Lagrangian splitting
$\d = \u + \up$, the intersection of any $N(\u)$-orbit and 
any $N(\up)$-orbit in $\Ld$ is a regular Poisson submanifold for
the Poisson structure $\Pi_{\u, \up}$ on $\Ld$.
\eth

To prove \thref{Ld-regular-2} 
we need to check that in the setting of 
\thref{Ld-regular-2} the condition of \prref{Ld-regular-1}
is satisfied. In fact, we prove a stronger statement.

\bpr{cond} If $\d$ is an even dimensional complex or real reductive Lie 
algebra and  $\lara$ is a nondegenerate symmetric 
invariant bilinear form on 
$\d$, then for all Lagrangian subalgebras $\l$ of $(\d, \lara)$,
$[\n(\l), \n(\l)] = (\n(\l))^\perp$.
\epr

The real case in \prref{cond} follows from the complex one.
Indeed, let $(\d, \lara)$ be a quadratic real reductive Lie algebra.
Then $(\d_{\Cset}, \lara_\Cset)$ is a quadratic complex reductive Lie algebra.
Let $\l$ be a Lagrangian subalgebra of $(\d, \lara)$, and let  
$\n(\l_\Cset)$ be the normalizer subalgebra of $\l_\Cset$ in $\d_\Cset$. Then
$\n(\l_\Cset) = \left(\n(\l)\right)_\Cset$. Assume 
the validity of \prref{cond} in the complex case.  We get
\begin{align*}
[\n(\l), \n(\l)] &= 
[\n(\l_\Cset), \n(\l_\Cset)] \cap \d =
\left(\n(\l_\Cset)\right)^\perp \cap \d \\
&= ( \left(\n(\l)\right)_\Cset )^\perp \cap \d =
(\n(\l))^\perp
\end{align*}
where $(.)^\perp$ denotes orthogonal
complements in $\d$ and $\d_\Cset$.
This proves the real case in \prref{cond}. 

To obtain the complex case in \prref{cond} we need 
the following result of Delorme \cite{delorme:manin-triples}.

\bth{de} [Delorme] Assume that $(\d, \lara)$ is an even dimensional 
reductive quadratic Lie algebra. For each Lagrangian subalgebra $\l$
of $(\d, \lara)$
the normalizer of the nilpotent radical $\n$ of $\l$ is 
a parabolic subalgebra $\p$ of $\d$. In addition, $\p$ has a
Levi subalgebra $\m$ whose derived subalgebra
$\bar{\m}=[\m,\m]$ decomposes as $\bar{\m} = \m_1 \oplus \m_2$
and for which there exists an isomorphism 
$\theta \colon \m_1 \rightarrow \m_2$. If $\z$ denotes the center 
of $\m$ then 
\begin{equation} 
\bar{\m}^{\theta} + \n \subset \l \subset 
\left( \bar{\m}^{\theta} \oplus \z \right) + \n
\label{dincl}
\end{equation}
where $\bar{\m}^{\theta} = \{ x + \theta(x) \mid x \in \m_1 \}
\subset \m_1 \oplus \m_2.$
\eth

\noindent
{\em{Proof of \prref{cond} in the complex case.}} 
Let $\l$ be a Lagrangian subalgebra of $(\d, \lara)$ as 
in \thref{de}. 

First we claim that $\n(\l) \subset \p$. Indeed, the normalizer 
of $\l$ lies inside the normalizer of the nilpotent radical $\n$ of $\l$ 
which is $\p$: if $y \in \n(\l)$, then for small $t$,
$\exp (t \ad_y)$ is an automorphism of $\l$ and thus 
of its nilpotent radical $\n$. Taking 
derivative at $t=0$, we get that $y$ normalizes $\n$. 

Next we will show that
\begin{equation}
\n(\l)= \left(\bar{\m}^{\theta} \oplus \z \right) + \n. 
\label{normll}
\end{equation}
The inclusion 
$\n(\l) \supset \left(\bar{\m}^{\theta} \oplus \z \right) + \n$ is clear
from \eqref{dincl}.
Define the subspace
\[
\bar{\m}^{-} = \{x - \theta(x) \mid x \in \m_1 \} \subset \bar{\m}.
\]
Under the adjoint action of $\bar{\m}^{\theta}$ we have 
the direct sum decomposition of $\bar{\m}^{\theta}$-modules
\[
\p = \bar{\m}^{\theta} \oplus \bar{\m}^{-} \oplus \z
\oplus \n.
\]
If $\n(\l) \neq \left(\bar{\m}^{\theta} \oplus \z \right) + \n$, 
then there exists a nonzero
$Y = y - \theta(y) \in \bar{\m}^{-}$ 
which belongs to $\n(\l)$. Since $\bar{\m}^{\theta}$ normalizes
$\bar{\m}^{-}$ we get that $\ad_Y(\bar{\m}^{\theta})=0$
and thus $\ad_y(\m_1)=0$. This is a contradiction
since $\m_1$ is semi-simple and $y \neq 0$. This
completes the proof of \eqref{normll}. 

Repeating the proof with $\n(\l)$ in the place of $\l$, 
which also satisfies \eqref{dincl} as shown above, we get
that $\n(\l)$ coincides with its normalizer.

Now \eqref{normll} 
implies $[\n(\l), \n(\l)] \subset \bar{\m}^{\theta} + \n.$
Because $\l$ is a Lagrangian subalgebra of $(\d, \lara)$, it is clear that
$\bar{\m}^{\theta}+ \n \subset \n(\l)^\perp$. Thus, 
$[\n(\l), \n(\l)] \subset \n(\l)^\perp$. By \reref{re-qqq}, 
$[\n(\l), \n(\l)] = \n(\l)^\perp$.
\qed
\sectionnew{Ranks of Poisson structures 
on the variety $\L(\g \oplus \g)$}
\label{Lgg}
\subsection{The quadratic Lie 
algebra $(\g \oplus \g, \lara)$}\lb{sec-gg}
Assume that $\g$ is a 
complex semi-simple Lie algebra and $\ll \, , \, \gg$ is a fixed
nondegenerate invariant symmetric
bilinear form whose restriction to a compact
real form of $\g$ is negative definite. 
Let $\d = \g \oplus \g$ be the direct
sum Lie algebra and let $\lara$ be the 
bilinear form on $\d$ given by
\begin{equation}
\lb{lara}
\la (x_1,x_2),(y_1,y_2)\ra =
\ll x_1,y_1\gg- \ll x_2,y_2\gg, \hspace{.2in} x_1, x_2, y_1, y_2
\in \g.
\end{equation}
In this section, we will study in more detail the Poisson structure $\Pill$
on $\Lgg$ defined by an arbitrary Lagrangian splitting $\gog = \l_1 + \l_2$.

A classification of Lagrangian subalgebras of $\g \oplus \g$ was first
obtained by Karolinsky \cite{karo:homog-diag}. It also follows from 
the more general results of Delorme \cite{delorme:manin-triples}, where
Lagrangian splittings of an 
arbitrary reductive quadratic Lie algebras were classified.
We will recall Delorme's
classification in $\S$\ref{sec-splitting-gg}. Let $G$ be the adjoint group 
of $\g$.
For $\l \in \L(\g \oplus \g)$
denote by $N(\l)$ the normalizer
subgroup of $\l$ in $G \times G$. Let $\gog = \l_1 + \l_2$
be an arbitrary Lagrangian splitting. By \thref{Ld-regular-2}
the intersection
of any $N(\l_1)$-orbit and any $N(\l_2)$-orbit in $\L(\g \oplus \g)$ is a 
regular
Poisson submanifold for the Poisson structure $\Pi_{\l_1, \l_1}$. Using 
results from 
\cite{LuY}, we will describe the $N(\l_1)$ and $N(\l_2)$-orbits in $\Lgg$
and will obtain an explicit formula for
the rank of $\Pill$ at an arbitrary $\l \in \Lgg$. 
\subsection{Lagrangian splittings of 
$(\g \oplus \g, \lara)$}\lb{sec-splitting-gg}
Fix a Cartan subalgebra $\h$ of $\g$ and a choice $\Delta^+$ of
positive roots in the set $\Delta$ of all roots for $(\g, \h)$. 
Let $\Gamma$ be the set of simple
roots in $\Delta^+$. For each $\alpha \in \Delta$, let $H_\alpha \in
\h$ be such that $\ll x, H_\alpha \gg =\alpha(x)$ for all $x \in \h$.
We will also fix a root vector $E_\alpha$ for each $\alpha \in \Delta$
such that $[E_\alpha, E_{-\alpha}] = H_\alpha$.
Following \cite{schiff:classcdyb, E-Lu2}, we define a {\it generalized 
Belavin--Drinfeld (gBD) triple} to be a triple $(S, T, d)$, where
$S$ and $T$ are subsets  of $\Gamma$ and $d: S \to T$ is a bijection
such that $\ll H_{d\alpha}, H_{d\beta} \gg = \ll H_\alpha, H_\beta\gg$
for all $\alpha \in S$. 

For a subset $S$ of $\Gamma$, let
$\Delta_S$ be the set of roots  
in the linear span of $S.$ Set
\[
\m_S = \h + {\sum}_{\alpha \in \Delta_S}^{} {\g}_\alpha, 
\hspace{.2in} \n_S = {\sum}_{\alpha\in {\De}^+ -\Delta_S}^{}
 {\g}_\alpha, \hs \n_{S}^{-} = \sum_{\alpha \in \Delta^+ - \Delta_S} 
\g_{-\alpha}
\]
and $\p_S = \m_S + \n_S$ and $\p_{S}^{-} = \m_S + \n_{S}^{-}.$
We set $\bar{\m}_S = [\m_S, \m_S]$ and
\[
\h_S = \h\cap\bar{\m}_S= \Span_\Cset \{H_\alpha: \alpha \in 
\Delta_S\}, \, \hs \, \z_S = \{x \in \h \mid \alpha(x) = 0, \,
\forall \alpha \in S\}.
\]
Then we have the decompositions 
$\h = \z_S +\h_S, \; \m_S = \z_S + \bar{\m}_S$ and
\[
\p_S = \z_S + \bar{\m}_S + \n_S, \, \hs \, \p_{S}^{-} = 
\z_S + \bar{\m}_S + \n_{S}^{-}.
\]
Recall that $G$ denotes the adjoint group of $\g$. 
The connected subgroups of $G$ with Lie algebras 
$\p_{S},$ $\p_{S}^{-}$, $\m_{S}$, $\n_S$ and $\n_{S}^{-}$ 
will be respectively denoted by
$P_{S}, P_{S}^{-}, M_{S}, N_S$ and $N_{S}^{-}$. 
Correspondingly we have the Levi decompositions $P_{S} = M_{S} N_{S}$,
$P_{S}^{-} = M_{S} N_{S}^{-}$.
Let $Z_S$ be the center of $M_S$, and let
$\chi_S: P_S \to M_S/Z_S$ be
the natural projection by first projecting to $M_S$ along $N_S$ and then to 
$M_S/Z_S$. We also denote by $\chi_S$  
the similar projection from $P_{S}^{-}$ to $M_S/Z_S$. 
 
For a generalized Belavin--Drinfeld triple $(S, T, d)$, 
let $\L_{\rm space}(\z_{S} \oplus \z_{T})$ be the set of all
Lagrangian subspaces of $\z_{S} \oplus \z_{T}$ with respect to the 
(nondegenerate) restriction of $\lara$ to $\z_{S} \oplus \z_{T}$. 
Let $\theta_d: \bar{\m}_S \to \bar{\m}_T$ be the unique Lie algebra 
isomorphism satisfying
\[
\theta_d(H_\alpha) = H_{d\al}, \hs 
\theta_d(E_\alpha) = E_{d\al}, \hs \forall \al \in S.
\]
For every $V \in \L_{\mathrm space}(\z_{S} \oplus \z_{T})$, define
\begin{align}\lb{lSTdV}
\l_{S, T, d, V} &= V + \{(x, \theta_d(x)) \mid x \in \bar{\m}_S\} 
+(\n_{S} \oplus \n_{T}) 
\subset \p_S \oplus \p_T,\\
\lb{lSTdV-prime}
\l_{S, T, d, V}^\prime &= V + \{(x, \theta_d(x)) \mid x \in \bar{\m}_S\} 
+(\n_{S} \oplus \n_{T}^-)
\subset \p_S \oplus \p_T^-,\\
\lb{lSTdV-2prime}
\l_{S, T, d, V}^{\prime\prime} &= 
V + \{(x, \theta_d(x)) \mid x \in \bar{\m}_S\}
+ (\n_{S}^- \oplus \n_{T})
\subset \p_S^- \oplus \p_T.
\end{align}
It is easy to see that $\l_{S, T, d, V}$, $\l^\prime_{S, T, d, V}$, and 
$\l^{\prime\prime}_{S, T, d, V}$ are all in $\Lgg$. 
The subalgebras $\l^\prime_{S, T, d, V}$ 
and $\l^{\prime\prime}_{S, T, d, V}$
are of course conjugate to ones of the type $\l_{S, T, d, V}$. Indeed, let
$W$ be the Weyl group of $(\g, \h)$, and let $w_0$ be the longest element in 
$W$. For $A \subset \Gamma$, let $W_A$ be the subgroup of $W$ generated by 
simple reflections with respect to roots in $A$, and
let $x_{A} = w_0 w_{0, A}$ where $w_{0, A}$ denotes the longest 
element of $W_A$. Then it is easy to see that
\begin{align}\lb{l-conj}
\l^\prime_{S, T, d, V} &= \Ad_{(e, \dot{x}_T)}^{-1} \l_{S, -w_0(T), \,x_T d,
\Ad_{(e, \dot{x}_T)}V},\\
\lb{l-conj-2}
\l^{\prime\prime}_{S, T, d, V} &= \Ad_{(\dot{x}_S, e)}^{-1} 
\l_{-w_0(S), T, \,dx_S^{-1}, \Ad_{(\dot{x}_S, e)}V},
\end{align}
where $\dot{x}_T$ and $\dot{x}_S$ are representatives in $G$ of $x_T$ and 
$x_S$ respectively. 
 
Denote also by $\theta_d$ the (unique) 
group isomorphism $M_S/Z_S \to M_T/Z_T$ induced by 
$\theta_d: \bar{\m}_S \to \bar{\m}_T$. 
Corresponding to the subalgebras in 
\eqref{lSTdV}-\eqref{lSTdV-2prime}, we define
\begin{align}\lb{RSTd}
R_{S, T, d} &= \{(p_1, p_2) \in P_S \times P_T \mid 
\theta_d(\chi_S(p_1)) = \chi_T(p_2)\} 
\subset P_S \times P_T, \\
\lb{RSTd-prime}
R_{S, T, d}^\prime &= \{(p_1, p_2) \in P_S \times P_T^- \mid 
\theta_d(\chi_S(p_1)) = \chi_T(p_2)\} 
\subset P_S \times P_T^-, \\
\lb{RSTd-2prime}
R_{S, T, d}^{\prime\prime} &= \{(p_1, p_2) \in P_S^- \times P_T \mid 
\theta_d(\chi_S(p_1)) = \chi_T(p_2)\} 
\subset P_S^- \times P_T.
\end{align}
One knows that $R_{S, T, d}$, $R_{S, T, d}^\prime$, and 
$R^{\prime\prime}_{S, T, d}$ are all connected \cite[Lemma 2.19]{E-Lu2}.
Corresponding to \eqref{l-conj} and \eqref{l-conj-2}, we have
\begin{align}\lb{R-conj-1}
R^\prime_{S, T, d} &= \Ad_{(e, \dot{x}_T)}^{-1} R_{S, -w_0(T), \,x_T d},\\
\lb{R-conj-2}
R^{\prime\prime}_{S, T, d} &= 
\Ad_{(\dot{x}_S, e)}^{-1} R_{-w_0(S), T, \,dx_S^{-1}}.
\end{align}
The Lie algebras of $R_{S, T, d}, R^\prime_{S, T, d}$, and 
$R^{\prime\prime}_{S, T, d}$ will
be denoted by $\r_{S, T, d}, \r^\prime_{S, T, d}$, and 
$\r^{\prime\prime}_{S, T, d}$
respectively. 

\bpr{GG-orbits}\cite{E-Lu2} 
Every $(G \times G)$-orbit in $\Lgg$ passes through an $\l_{S, T, d, V}$ for
a unique generalized Belavin--Drinfeld triple $(S, T, d)$ and a unique
$V \in \L_{\rm space}(\z_{S} \oplus \z_{T})$. The normalizer subgroup of
$\l_{S, T, d,V}$ in $G \times G$ is $R_{S, T, d}$.
\epr

\bde{BD-system} 
For generalized Belavin-Drinfeld triples 
$(S_i, T_i, d_i)$, $i = 1, 2$, let
\[
S_2^{d_1^{-1}d_2} = \{\alpha \in S_2 \mid \;
(d_1^{-1}d_2)^n \alpha \; 
\mbox{is defined and is in} \; S_2 \; \mbox{for} \; n = 1, 2, \cdots\}.
\]
A {\bf generalized Belavin--Drinfeld system} is a pair of quadruples
$(S_1, T_1, d_1, V_1)$ and $(S_2, T_2, d_2, V_2)$, where for $i = 1, 2$, 
$(S_i, T_i, d_i)$ is a generalized Belavin--Drinfeld triple and 
$V_i \in \L_{\rm space}(\z_{S_i} \oplus \z_{T_i})$,
such that

1) $S_2^{d_1^{-1}d_2} = \emptyset$;

2) $\h_1 \cap \h_2 = \{0\}$, where 
$\h_i = V_i + \{(x, \theta_d(x)) \mid x \in \h_{S_i}\} 
\subset \h \oplus \h$ for $i = 1, 2$.
\ede

\bth{delorme-splitting}\cite[Delorme]{delorme:manin-triples}
Every Lagrangian splitting of $\gog$ is conjugate by 
an element in $G \times G$
to one of the form $\gog = \l_1 + \l_2$, where 
\begin{equation}\lb{l1l2}
\l_1 = \l_{S_1, T_1, d_1, V_1}^\prime 
\hspace{.2in} \mbox{and} \hspace{.2in}
\l_2 = \l_{S_2, T_2, d_2, V_2}^{\prime\prime}
\end{equation}
for a  generalized Belavin--Drinfeld system
$(S_1, T_1, d_1, V_1)$, $(S_2, T_2, d_2, V_2)$.
\eth

Let $\gog = \l_1 + \l_2$ be a Lagrangian splitting with $\l_1$ and $\l_2$ 
given in \eqref{l1l2}.
By \thref{Ld-regular-2}, any non-empty intersection of an $N(\l_1)$
and an $N(\l_2)$-orbit in $\Lgg$ is a regular Poisson subvariety for
the Poisson structure $\Pi_{\l_1, \l_2}$.
A classification of $N(\l_1)$ and $N(\l_2)$-orbits will be given in 
$\S$ \ref{Nl-orbits}. We now prove that every non-empty intersection
of an $N(\l_1)$-orbit and an $N(\l_2)$-orbit in $\Lgg$ is smooth and 
irreducible. Let $H$ be the connected subgroup of $G$ with Lie algebra $\h$.

\bpr{NN-connected}
For a Lagrangian splitting $\gog = \l_1 + \l_2$ with $\l_1 $ and $\l_2$ 
given by \eqref{l1l2}, $N(\l_1) \cap N(\l_2)$ is a subtorus of $H \times H$ 
of dimension $\dim \z_{S_1} + \dim \z_{S_2}$. In particular, 
$N(\l_1) \cap N(\l_2)$ is connected.
\epr

\begin{proof}
It follows from \prref{GG-orbits} that 
\begin{equation}\lb{NN}
N(\l_1) = R_{S_1, T_1, d_1}^\prime \hspace{.2in} \mbox{and} \hspace{.2in}
N(\l_2) = R_{S_2, T_2, d_2}^{\prime\prime}.
\end{equation}
For notational simplicity, let $x_1 = x_{T_1}, x_2 = x_{S_2}$. 
By \eqref{R-conj-1} and \eqref{R-conj-2}, 
\[
N(\l_1) \cap N(\l_2) =\left( \Ad_{(e, \dot{x}_{1})}^{-1}
R_{S_1, -w_0(T_1), \,x_{1}d_1}\right) \cap \left( \Ad_{(\dot{x}_{2}, e)}^{-1}
R_{-w_0(S_2), \, T_2, \,d_2x_{2}^{-1}}\right).
\]
Since $x_2^{-1} \in W^{-w_0(S_2)}$ and
$x_1 \in {}^{-w_0(T_1)}\!W$, 
we can use \cite[Theorem 2.5]{LuY} to determine $N(\l_1) \cap N(\l_2)$. 
Set $N = N_\emptyset$.
Since
$S_2^{d_1^{-1}d_2} = \emptyset$, \cite[Theorem 2.5]{LuY} implies that
$N(\l_1) \cap N(\l_2) \subset B \times B$ and
$N(\l_1) \cap N(\l_2)= (N(\l_1) \cap N(\l_2))^{{\rm red}} 
(N(\l_1) \cap N(\l_2))^{{\rm uni}},$
where 
\begin{align*}
(N(\l_1) \cap N(\l_2))^{{\rm red}} &= 
N(\l_1) \cap N(\l_2) \cap (H \times H), \\
(N(\l_1) \cap N(\l_2))^{{\rm uni}} &= 
N(\l_1) \cap N(\l_2) \cap (N \times N).
\end{align*}
Moreover, \cite[Theorem 2.5]{LuY} also tells us that
\[
(N(\l_1) \cap N(\l_2))^{{\rm uni}} \cong 
(N \cap \Ad_{\dot{x}_2}^{-1} (N_{-w_0(S_2)}) ) 
\times (N \cap \Ad_{\dot{x}_1}^{-1}(N_{-w_0(T_1)})).
\]
It is easy to see that $(N \cap \Ad_{\dot{x}_2}^{-1} (N_{-w_0(S_2)}) ) 
\times (N \cap \Ad_{\dot{x}_1}^{-1}(N_{-w_0(T_1)}))$ is the trivial group. 
Thus $N(\l_1) \cap N(\l_2) = N(\l_1) \cap N(\l_2) \cap (H \times H)$ consists 
of all $(h_1, h_2) \in H \times H$ such that
\[
\begin{cases}
& \theta_{d_1} \chi_{S_1}(h_1) = 
\chi_{T_1}(h_2), \\
& \theta_{d_2} \chi_{S_2} (h_1) = \chi_{T_2}(h_2),
\end{cases}
\]
which are equivalent to 
\begin{equation}\lb{hh}
h_1^\alpha = h_2^{d_1\alpha}, \;\;  \forall \alpha \in S_1 \; \; \; 
\mbox{and} \; \; \; 
h_1^\beta = h_2^{d_2\beta}, \; \;\forall \beta \in S_2.
\end{equation}
Let $\Gamma = \{\alpha_1, \, \alpha_2, \, \cdots, \, \alpha_r\}$ be the set 
of simple roots of $\g$.
Since $G$ is the adjoint group of $\g$, 
we can identify $H \times H$ with the torus $(\Cset^\times)^{2r}$ by the map
\[
H \times H \longrightarrow (\Cset^\times)^{2r}: \; \; \; (h_1, h_2) 
\longmapsto (h_1^{\alpha_1}, \;
h_1^{\alpha_2}, \;\cdots, \;h_1^{\alpha_r}, \;
h_2^{\alpha_1}, \;h_2^{\alpha_2}, \;\cdots, \;h_2^{\alpha_r}).
\]
The conditions in \eqref{hh} imply that the 
coordinates $h_1^\alpha$ of $h_1$ for
$\alpha \in S_1 \cup S_2$ are expressed in terms of coordinates of
$h_2$, and we have the extra conditions  
\begin{equation}\lb{hh-extra}
h_2^{d_1\alpha} = h_2^{d_2\alpha}, \; \; \; \alpha \in S_1 \cap S_2
\end{equation}
for the coordinates of $h_2$. To understand the conditions in 
\eqref{hh-extra}, recall that
$S_2^{d_1^{-1}d_2} = \emptyset$. Thus for every $\alpha \in S_1 \cap S_2$, 
there is a unique integer $n \geq 1$ and unique elements $\alpha^{(0)} 
=\alpha, \, \alpha^{(1)}, \, \alpha^{(2)}, \, \cdots, \, 
\alpha^{(n-1)} \in S_1 \cap S_2$
such that 
\[
d_2 \alpha^{(0)} = d_1\alpha^{(1)}, \; \;  
d_2 \alpha^{(1)} = d_1\alpha^{(2)}, \; \;  \cdots \; \; 
d_2 \alpha^{(n-2)} = d_1\alpha^{(n-1)}
\]
and either $d_2 \alpha^{(n-1)} \notin T_1$ or $d_2\alpha^{(n-1)} =
d_1\alpha^{(n)} $ for some $\alpha^{(n)} \in S_1$ but $\alpha^{(n)} \notin 
S_2$.
Then the conditions in \eqref{hh-extra} are equivalent to 
\[
h_2^{d_1\alpha} = h_2^{d_2\alpha} = 
h_2^{d_2\alpha^{(1)}} = \cdots = h_2^{d_2\alpha^{(n-1)}}.
\]
Since $d_2 \alpha^{(n-1)} \notin d_1(S_1 \cap S_2)$, we see that the 
conditions \eqref{hh-extra} express $h_2^{d_1\alpha}$ for every 
$\alpha \in S_1 \cap S_2$
in terms of $h_2^\beta$ for some $\beta \notin d_1(S_1 \cap S_2)$. We 
conclude that the set of $(h_1, h_2) \in H \times H$ satisfying \eqref{hh} 
is a subtorus of $H \times H$ with dimension equal to
\[
2 \dim H - |S_1 \cup S_2| - |S_1 \cap S_2| = 2\dim H -|S_1|-|S_2| = 
\dim \z_{S_1} + \dim \z_{S_2}.
\] 
\end{proof}

\bco{co-NN-irred}
For any Lagrangian splitting $\gog = \l_1 + \l_2$, all $N(\l_1)$-orbits and
$N(\l_2)$-orbits in $\L(\gog)$ intersection transversally, and every such 
non-empty intersection is smooth and irreducible.
\eco

\begin{proof} Clearly $\n(\l_1) + \n(\l_2) = \gog$, where $\n(\l_i)$ is
the Lie algebra of $N(\l_i)$ for $i = 1, 2$. Since $N(\l_1) \cap N(\l_2)$
is connected, \coref{co-NN-irred} follows from \cite[Corollary 1.5]{Ri}.
\end{proof}
\subsection{$N(\l_1)$ and $N(\l_2)$-orbits in $\Lgg$}\lb{Nl-orbits}
Assume that $\gog=\l_1 + \l_2$ is a Lagrangian splitting with $\l_1$ and 
$\l_2$ given in \eqref{l1l2}. We now use results in \cite{LuY} to describe 
$N(\l_1)$ and $N(\l_2)$-orbits in $\Lgg$.

For $A \subset \Gamma$, let $W^A$ and ${}^A\!W$ be respectively the 
sets of minimal length representatives 
in the cosets in $W/W_A$ and $W_A \backslash W$.
For each $w \in W$, we fix a representative $\dot{w}$ of $w$ in 
the normalizer of $H$ in $G$.

\bpr{N-orbits}
1) Every $N(\l_1) =R^\prime_{S_1, T_1, d_1}$-orbit in $\Lgg$ is of the form
\[
R^\prime_{S_1, T_1, d_1} \Ad_{(\dot{v}_1, \dot{v}_2m_2)}\l_{S, T, d, V}
\]
for a unique generalized Belavin-Drinfeld triple $(S, T, d)$, a unique
$V \in {{\mathcal L}}_{{\rm space}}(\z_{S} \oplus \z_T)$, a unique
pair $(v_1, v_2) \in W^S \times {}^{T_1}\!W$, and some $m_2 \in M_{T(v_1, v_2)}$ 
with
\[
T(v_1, v_2) = \{\alpha \in T \mid (v_2^{-1} d_1 v_1 d^{-1})^n \alpha \; 
\mbox{is defined and is in} \; T\; \mbox{for} \; n = 1, 2, \cdots\}.
\]

2) Every $N(\l_2)=R^{\prime\prime}_{S_2, T_2, d_2}$-orbit in $\Lgg$  
is of the form
\[
R^{\prime\prime}_{S_2, T_2, d_2}
\Ad_{(\dot{w}_1 m_1, \dot{w}_2)} \l_{S, T, d, V}
\]
for a unique generalized Belavin-Drinfeld triple $(S, T, d)$, a unique
$V \in {{\mathcal L}}_{{\rm space}}(\z_{S} \oplus \z_T)$, a unique
pair $(w_1, w_2) \in {}^{S_2}\!W \times W^T$, and some $m_1 \in M_{S(w_1, w_2)}$ 
with
\[
S(w_1, w_2) = \{\alpha \in S \mid  (w_1^{-1} d_2^{-1} w_2 d)^n \alpha \; 
\mbox{is defined and is in} \; S\; \mbox{for} \; n = 1, 2, \cdots\}.
\]
\epr


\begin{proof} In \cite{LuY} we gave a 
description of the $(R_{S_1, T_1, d_1},
R_{S, T, d})$-double cosets in $G \times G$. 
The first part follows from \eqref{R-conj-1}, 
\cite[Theorem 2.2]{LuY}, and the
fact that ${}^{-w_0(A)}\!W = x_A  {}^{A}\!W$, 
$\forall A \subset \Gamma$, see \S \ref{sec-splitting-gg} 
for the definition of $x_A$.

Let $\sigma: G \times G \to G \times G: (g_1, g_2) \mapsto (g_2, g_1)$. 
Then using first part and the facts that
\[
R^{\prime\prime}_{S_2, T_2, d_2} = 
\sigma \left( R^\prime_{T_2, S_2, d_2^{-1}} \right)
\hspace{.2in} \mbox{and} \hspace{.2in} 
R_{S, T, d} = \sigma ( R_{T, S, d^{-1}} )
\]
we get part 2).
\end{proof}

Let $\O_1$ be an $N(\l_1) = R^\prime_{S_1, T_1, d_1}$-orbit 
and $\O_2$ an $N(\l_2) = R^{\prime\prime}_{S_2, T_2, d_2}$-orbit in
$\Lgg$. By \prref{N-orbits}, we can assume that
\begin{align}\lb{O1}
\O_1 &= R^\prime_{S_1, T_1, d_1} 
\Ad_{(\dot{v}_1, \dot{v}_2 m_2)} \l_{S, T, d, V},\\
\lb{O2}\O_2 & = R^{\prime\prime}_{S_2, T_2, d_2} 
\Ad_{(\dot{w}_1 m_1, \dot{w}_2)} \l_{S, T, d, V},
\end{align}
where $(S, T, d, V)$, $v_1, v_2, w_1, w_2$ and $m_1$ and $m_2$ are as in 
\prref{N-orbits}. Let
\begin{align}\lb{S1vv}
S_1(v_1, v_2)& = d_1^{-1} v_2 T(v_1, v_2) = 
v_1 d^{-1} T(v_1, v_2) \subset S_1\\
S_2(w_1, w_2) &=w_1S(w_1, w_2) = 
d_2^{-1} w_2 d S(w_1, w_2) \subset S_2.
\end{align}
In other words, $S_1(v_1, v_2)$ is the largest 
subset of $S_1$ that is invariant under
the partial map $v_1 d^{-1} v_2^{-1} d_1: \Gamma \to \Gamma$, and
$S_2(w_1, w_2)$ is the largest subset of $S_2$ 
that is invariant under the partial map
$w_1 d^{-1} w_2^{-1} d_2: \Gamma \to \Gamma$.

In order to compute the 
rank of $\Pi_{\l_1, \l_2}$ at an $\l \in \O_1 \cap \O_2$ using \leref{formula-NN}, we need to
compute the dimensions of various intersections of subalgebras in $\n(\l_1), \n(\l_2)$ and $\n(\l)$.
Such intersections can be described using the following 
\prref{from-LuY} which is derived from \cite[Theorem 2.5]{LuY}.

Recall that for 
$\q \subset \gog, \q^\perp = \{ x \in \gog \mid \la x, y \ra = 0, \; \forall y \in \q\}$.
Clearly
\begin{align*}
\r^\perp_{S, T, d} &= \n_S \oplus \n_T + \{(x, \theta_d(x))\mid x \in \bar{\m}_S\},\\
\r^{\prime, \perp}_{S_1, T_1, d_1}  &= \n_{S_1} \oplus \n_{T_1}^- + 
\{(x, \theta_{d_1}(x)) \mid x \in \bar{\m}_{S_1}\},\\
\r^{\prime\prime, \perp}_{S_2, T_2, d_2} &=
\n_{S_2}^{-} \oplus \n_{T_2} + \{(x, \theta_{d_2}(x)) \mid x \in \bar{\m}_{S_2}\}.
\end{align*}
\bnota{aaa}
Let $\calS$ (reps. $\calS^\prime$, $\calS^{\prime\prime}$) be the set of all subspaces of 
$\r_{S, T, d}$ (resp. $\r^\prime_{S_1, T_1, d_1}$, $\r^{\prime\prime}_{S_2, T_2, d_2}$)
that contain $\r^\perp_{S, T, d}$ (resp. $\r^{\prime, \perp}_{S_1, T_1, d_1}$,
$\r^{\prime\prime, \perp}_{S_2, T_2, d_2}$).
For $\a  \in \calS$, let $V_\a = \a \cap (\z_S \oplus \z_T)$ and
\begin{align*}
X_\a^\prime &= V_\a + \{(x, \theta_d(x)) \mid x \in \z_{d^{-1}T(v_1, v_2)} \cap \h_S\},\\
X_\a^{\prime\prime} &= V_\a + \{(x, \theta_d(x)) \mid x \in \z_{S(w_1, w_2)} \cap \h_S\}.
\end{align*}
For $\a^\prime \in \calS^\prime$ and $\a^{\prime\prime} \in \calS^{\prime\prime}$, let
\begin{align*}
V_{\a^\prime}& = \a^\prime \cap (\z_{S_1} \oplus \z_{T_1}) \hs \mbox{and} \hs
Y_{\a^\prime} = V_{\a^\prime} + \{(x, \theta_{d_1}(x)) \mid x \in 
\z_{S_1(v_1, v_2)} \cap \h_{S_1}\},\\
V_{\a^{\prime\prime}}& = \a^{\prime\prime} \cap (\z_{S_2} \oplus \z_{T_2}) \hs 
\mbox{and} \hs
Y_{\a^{\prime\prime}} = V_{\a^{\prime\prime}} + \{(x, \theta_{d_2}(x)) 
\mid x \in \z_{S_2(w_1, w_2)} \cap \h_{S_2}\}.
\end{align*}
We also set
\begin{equation}\lb{f1}
\f_1 = \r^{\prime, \perp}_{S_1, T_1, d_1} \cap \Ad_{(\dot{v}_1, \dot{v}_2m_2)} 
\r^\perp_{S, T, d} \cap (\bar{\m}_{S_1(v_1, v_2)} \oplus 
\bar{\m}_{d_1S_1(v_1, v_2)}).
\end{equation}
\enota

\bpr{from-LuY}
For any $\a^\prime \in \calS_1$ and $\a \in \calS$, one has the direct sum
\[
\a^\prime \cap \Ad_{(\dot{v}_1, \dot{v}_2 m_2)} \a = (\a^\prime \cap 
\Ad_{(\dot{v}_1, \dot{v}_2 m_2)} \a)^{\rm red}
+ (\a^\prime \cap \Ad_{(\dot{v}_1, \dot{v}_2 m_2)} \a)^{\rm nil},
\]
where 
\begin{align*}
(\a^\prime \cap \Ad_{(\dot{v}_1, \dot{v}_2 m_2)} \a)^{\rm red} &=
\a^\prime \cap \Ad_{(\dot{v}_1, \dot{v}_2 m_2)} \a \cap (\m_{S_1(v_1, v_2)} \oplus 
\m_{d_1S_1(v_1, v_2)})\\
& = Y_{\a^\prime} \cap \Ad_{(\dot{v}_1, \dot{v}_2)} X_\a^\prime + \f_1 \hs (\mbox{direct sum}),
\end{align*}
and $(\a^\prime \cap \Ad_{(\dot{v}_1, \dot{v}_2 m_2)} \a)^{\rm nil} =
\a^\prime \cap \Ad_{(\dot{v}_1, \dot{v}_2 m_2)} \a \cap (\n_{S_1(v_1, v_2)} \oplus 
\n_{d_1S_1(v_1,v_2)}^-)$. The dimension of the latter is equal to
$l(v_2) + \dim(\n \cap \Ad_{\dot{v}_1}(\n_S))$.
\epr

\begin{proof}
Using \eqref{l-conj} and \cite[Theorem 2.5]{LuY}, one can see that
\[
\a^\prime \cap \Ad_{(\dot{v}_1, \dot{v}_2 m_2)} \a \subset
\p_{S_1(v_1, v_2)} \oplus \p_{d_1S_1(v_1, v_2)}^-
\]
and that 
\begin{align*}
\a^\prime \cap \Ad_{(\dot{v}_1, \dot{v}_2 m_2)} \a &= 
(\a^\prime \cap \Ad_{(\dot{v}_1, \dot{v}_2 m_2)} \a)^{\rm red}
+ (\a^\prime \cap \Ad_{(\dot{v}_1, \dot{v}_2 m_2)} \a)^{\rm nil}\\
&=\a^\prime \cap \Ad_{(\dot{v}_1, \dot{v}_2 m_2)} \a \cap (\m_{S_1(v_1, v_2)} \oplus 
\m_{d_1S_1(v_1, v_2)}) \\
& \; \; +\a^\prime \cap \Ad_{(\dot{v}_1, \dot{v}_2 m_2)} \a \cap (\n_{S_1(v_1, v_2)} \oplus 
\n_{d_1S_1(v_1,v_2)}^-).
\end{align*}
The dimension formula for  $(\a^\prime \cap \Ad_{(\dot{v}_1, \dot{v}_2 m_2)} \a)^{\rm nil}$
also follows from \cite[Theorem 2.5]{LuY}. It now remains to show that 
$(\a^\prime \cap \Ad_{(\dot{v}_1, \dot{v}_2 m_2)} \a)^{\rm red}
=Y_{\a^\prime} \cap \Ad_{(\dot{v}_1, \dot{v}_2)} X_\a^\prime + \f_1$ as a direct sum.
Clearly, $Y_{\a^\prime} \cap \Ad_{(\dot{v}_1, \dot{v}_2)} X_\a^\prime$ and $\f_1$ intersect 
trivially, and their sum is contained in 
$(\a^\prime \cap \Ad_{(\dot{v}_1, \dot{v}_2 m_2)} \a)^{\rm red}$.

Suppose that $(x, y) \in (\a^\prime \cap \Ad_{(\dot{v}_1, \dot{v}_2 m_2)} \a)^{\rm red}$.
Write $x = x_1 + x_2$ and $y = y_1 + y_2$, where
$x_1 \in \z_{S_1(v_1, v_2)}, x_2 \in \bar{\m}_{S_1(v_1, v_2)},
y_1 \in \z_{d_1S_1(v_1, v_2)}$ and $
y_1 \in \bar{\m}_{d_1S_1(v_1, v_2)}$. Because
$(x, y) \in \a^\prime \subset \r^\prime_{S_1, T_1, d_1}$, we have 
$\theta_{d_1} \chi_{S_1}(x_1) + \theta_{d_1}(x_2) = \chi_{T_1}(y_1) + y_2$. 
It follows from the direct sum decomposition
$\z_{d_1S_1(v_1, v_2)} = \z_{T_1} + \z_{d_1S_1(v_1, v_2)} \cap \h_{T_1}$ that
$\chi_{T_1}(y_1) \in \z_{S_1(v_1, v_2)} \cap \h_{T_1}$. Similarly,
$\chi_{S_1}(x_1) \in \z_{S_1(v_1, v_2)} \cap \h_{S_1}$. 
Thus $\theta_{d_1}\chi_{S_1}(x_1) = \chi_{T_1}(y_1)$ and $\theta_{d_1}(x_2) = y_2$.
Hence $(x_2, y_2) \in \r^{\prime,\perp}_{S_1, T_1, d_1} \subset \a^\prime$, and therefore
$(x_1, y_1) \in (\z_{S_1(v_1, v_2)} \oplus \z_{d_1S_1(v_1, v_2)}) \cap \a^\prime
=Y_{\a^\prime}$.
In the same way one shows that $(x_1, y_1) \in \Ad_{(\dot{v}_1, \dot{v}_2)} X_{\a}^\prime$
and $(x_2, y_2) \in  \Ad_{(\dot{v}_1, \dot{v}_2m_2)} \r^{\perp}_{S, T, d}$.
Thus $(x_1, y_1) \in Y_{\a^\prime} \cap \Ad_{(\dot{v}_1, \dot{v}_2)} X_\a^\prime$
and $(x_2, y_2) \in \f_1$.
\end{proof}

\bco{a-b}
For any $\a, \b \in \calS$, $\a^\prime, \b^\prime \in \calS^\prime$, and $\a^{\prime\prime}, 
\b^{\prime\prime} \in \calS^{\prime\prime}$,
\begin{align}\lb{1}
\dim &(\a^\prime \cap \Ad_{(\dot{v}_1, \dot{v}_2 m_2)} \a) - 
\dim(\b^\prime \cap \Ad_{(\dot{v}_1, \dot{v}_2 m_2)} \b)\\ \nonumber
&=\dim(Y_{\a^\prime} \cap \Ad_{(\dot{v}_1, \dot{v}_2)} X_\a^\prime) - 
\dim(Y_{\b^\prime} \cap \Ad_{(\dot{v}_1, \dot{v}_2)} X_\b^\prime),\\
\lb{2}
\dim &(\a^{\prime\prime} \cap \Ad_{(\dot{w}_1 m_1, \dot{w}_2)} \a) - 
\dim(\b^{\prime\prime} \cap \Ad_{(\dot{w}_1 m_1, \dot{w}_2)} \b)\\ \nonumber
&=\dim(Y_{\a^{\prime\prime}} \cap \Ad_{(\dot{w}_1, \dot{w}_2)} X_\a^{\prime\prime}) - 
\dim(Y_{\b^{\prime\prime}} \cap \Ad_{(\dot{w}_1, \dot{w}_2)} X_\b^{\prime\prime}).
\end{align}
\eco

\begin{proof}
Let $\f_1$ be as in\eqref{f1}.  
We know from \prref{from-LuY} that 
\[
(\a^\prime \cap \Ad_{(\dot{v}_1, \dot{v}_2 m_2)} \a)^{\rm red} = 
Y_{\a^\prime} \cap \Ad_{(\dot{v}_1, \dot{v}_2)} X_\a^\prime + \f_1
\]
is a direct sum. 
Replacing $\a^\prime$ by $\b^\prime$ and $\a$ by $\b$, we get
\begin{align*}
\dim &(\a^\prime \cap \Ad_{(\dot{v}_1, \dot{v}_2 m_2)} \a) - 
\dim(\b^\prime \cap \Ad_{(\dot{v}_1, \dot{v}_2 m_2)} \b)\\
& \; \; \; =
\dim (\a^\prime \cap \Ad_{(\dot{v}_1, \dot{v}_2 m_2)} \a)^{\rm red} - 
\dim(\b^\prime \cap \Ad_{(\dot{v}_1, \dot{v}_2 m_2)} \b)^{\rm red}\\
& \; \; \; =
 \dim(Y_{\a^\prime} \cap \Ad_{(\dot{v}_1, \dot{v}_2)} X_\a^\prime) - 
\dim(Y_{\b^\prime} \cap \Ad_{(\dot{v}_1, \dot{v}_2)} X_\b^\prime).
\end{align*}
\eqref{2} is proved by using \eqref{1} and the map $\sigma: \gog \to \gog:
(x, y) \mapsto (y, x)$.
\end{proof}

\subsection{The rank of the Poisson structure $\Pi_{\l_1,\l_2}$ on $\Lgg$}
\label{sec_rank-Lgg}


\bth{rank-main}
Let $\gog = \l_1 + \l_2$ be a Lagrangian splitting as in \eqref{l1l2}
and let $\O_1$ and $\O_2$
be respectively an $N(\l_1)$ and an $N(\l_2)$-orbit in
$\Lgg$ as in \eqref{O1} and \eqref{O2}. 
Then the corank of $\Pi_{\l_1, \l_2}$ in $\O_1 \cap \O_2$ is equal to
\begin{align*}
\dim \z_{S_1} + \dim \z_{S_2} + \dim \z_S 
&-\dim(Y_1 \cap \Ad_{(\dot{v}_1, \dot{v}_2)} X_1) +
\dim (Z_1 \cap \Ad_{(\dot{v}_1, \dot{v}_2)} X_1)\\
& -\dim(Y_2 \cap \Ad_{(\dot{w}_1, \dot{w}_2)} X_2)+
\dim (Z_2 \cap 
\Ad_{(\dot{w}_1, \dot{w}_2)} \widetilde{X}),
\end{align*}
where
\begin{align*}
X_1 &= (\z_{d^{-1}T(v_1, v_2)} \oplus \z_{T(v_1, v_2)}) \cap r_{S, T, d}\\
&=\z_{S} \oplus \z_T \!+\! \{(x, \theta_d(x)) \mid x \in \z_{d^{-1}T(v_1, v_2)}\cap \h_S\}\\
X_2 &=(\z_{S(w_1, w_2)} \oplus \z_{dS(w_1, w_2)}) \cap r_{S, T, d}\\
&=\z_{S} \oplus \z_T \!+\! \{(x, \theta_d(x)) \mid x \in \z_{S(w_1, w_2)}\cap \h_S\}\\
Y_1 &=(\z_{S_1(v_1, v_2)} \oplus \z_{d_1S_1(v_1, v_2)}) \cap \r^\prime_{S_1, T_1, d_1}\\
&=\z_{S_1} \oplus \z_{T_1} + \{(x, \theta_{d_1}(x)) \mid x \in \z_{S_1(v_1, v_2)} \cap \h_{S_1}\} \\
Y_2 &=(\z_{S_2(w_1, w_2)} \oplus \z_{d_2S_2(w_1, w_2)}) \cap  
\r^{\prime\prime}_{S_2, T_2, d_2}\\
&=\z_{S_2} \oplus \z_{T_2} + \{(x, \theta_{d_2}(x)) \mid x \in \z_{S_2(w_1, w_2)} \cap \h_{S_2} \} \\
Z_1&= (\z_{S_1(v_1, v_2)} \oplus \z_{d_1S_1(v_1, v_2)}) \cap \l^\prime_{S_1, T_1, d_1, V_1}\\
&=V_1 + \{(x, \theta_{d_1}(x)) \mid x \in \z_{S_1(v_1, v_2)} \cap \h_{S_1}\} \\
Z_2 &=(\z_{S_2(w_1, w_2)} \oplus \z_{d_2S_2(w_1, w_2)}) \cap  
\l^{\prime\prime}_{S_2, T_2, d_2, V_2}\\
&=V_2 + \{(x, \theta_{d_2}(x)) \mid x \in \z_{S_2(w_1, w_2)} \cap \h_{S_2} \},
\end{align*}
and
$\tilde{X} = p(X_1\cap \Ad_{(\dot{v}_1, \dot{v}_2)}^{-1}Z_1) + 
\{(x, \theta_d(x)) \mid x \in \z_{S(w_1, w_2)}\cap \h_S\}$ with
$p: \hoh \to \z_S \oplus \z_T$ being the projection along $\h_S \oplus \h_T$.
\eth

\begin{proof} Let $\l \in \O_1 \cap \O_2$ be given by 
\begin{equation}\lb{l}
\l = \Ad_{(r_1^\prime, r_2^\prime)(\dot{v}_1, \dot{v}_2m_2)}\l_{S, T, d, V}
=\Ad_{(r_1^{\prime\prime}, r_2^{\prime\prime})(\dot{w}_1m_1, \dot{w}_2)}\l_{S, T, d, V},
\end{equation}
where $(r_1^\prime, r_2^\prime) \in N(\l_1) = R^\prime_{S_1, T_1, d_1}$ and
$(r_1^{\prime\prime}, r_2^{\prime\prime}) \in N(\l_2) = R^{\prime\prime}_{S_2, T_2, d_2}$.
The formula for the corank of $\Pi_{\l_1, \l_2}$ at $\l \in 
\O_1 \cap \O_2$  involves the Drinfeld subalgebra 
${\mathcal{T}}(\l)$ of
$\gog$ defined by ${\mathcal{T}}(\l) = \n(\l)^\perp + \l_1 \cap \n(\l)$,
cf. \prref{ld},
where again $\n(\l)$ is the normalizer subalgebra of $\l$ in $\gog$.
We first compute ${{\mathcal T}}(\l)$. 
Since 
\[
\Ad_{(r_1^\prime, r_2^\prime) 
(\dot{v}_1, \dot{v}_2m_2)}\r^\perp_{S, T, d}=
\n(\l)^\perp \subset {{\mathcal T}}(\l) \subset \n(\l)=
\Ad_{(r_1^\prime, r_2^\prime) 
(\dot{v}_1, \dot{v}_2m_2)}\r_{S, T, d},
\]
we know that ${{\mathcal T}}(\l)=\Ad_{(r_1^\prime, r_2^\prime) 
(\dot{v}_1, \dot{v}_2m_2)}\l_{S, T, d, \widetilde{V}}$ 
for some $\widetilde{V} \in \L(\z_S \oplus \z_T)$. On the other hand,
\[
{{\mathcal T}}(\l) =\Ad_{(r_1^\prime, r_2^\prime) 
(\dot{v}_1, \dot{v}_2m_2)} \left( \r^\perp_{S, T, d} + r_{S, T, d} \cap 
\Ad_{(\dot{v}_1, \dot{v}_2 m_2)}^{-1} \l_1\right).
\]
Thus $\l_{S, T, d, \widetilde{V}} = 
\r^\perp_{S, T, d} + r_{S, T, d} \cap 
\Ad_{(\dot{v}_1, \dot{v}_2 m_2)}^{-1} \l_1=
X_1 \cap \Ad_{(\dot{v}_1, \dot{v}_2)}^{-1}Z_1 + \r^\perp_{S, T, d}$,
where the second identity comes from \prref{from-LuY}. Hence
\[
\widetilde{V}= p(X_1 \cap \Ad_{(\dot{v}_1, \dot{v}_2)}^{-1}Z_1).
\]
Now by \leref{formula-NN},
the corank
of $\Pi_{\l_1, \l_2}$ in $\O_1 \cap \O_2$
at the Lagrangian subalgebra $\l$ given by \eqref{l} is
\begin{align*}
&{\rm Corank}_{\Pi_{\l_1, \l_2}}(\l) = \dim \z_{S_1} + 
\dim \z_{S_2} + \dim \z_S \\
& \; \; \; \; -\dim(\r^\prime_{S_1, T_1, d_1} \cap 
\Ad_{(\dot{v}_1, \dot{v}_2 m_2)} \r_{S, T, d})
+ \dim(\l^\prime_{S_1, T_1, d_1, V_1} \cap 
\Ad_{(\dot{v}_1, \dot{v}_2 m_2)} \r_{S, T, d})\\
&\; \; \; \; -\dim(\r^{\prime\prime}_{S_2, T_2, d_2} \cap 
\Ad_{(\dot{w}_1m_1, \dot{w}_2)} \r_{S, T, d})
+\dim(\l^{\prime\prime}_{S_2, T_2, d_2} \cap 
\Ad_{(\dot{w}_1m_1, \dot{w}_2)} \l_{S, T, d, \widetilde{V}}).
\end{align*}
Applying \coref{a-b}, we get the 
desired formula for the corank of $\Pi_{\l_1, \l_2}$ in $\O_1 \cap \O_2$.
This completes the proof of \thref{rank-main}.
\end{proof}

\bex{Belavin-Drinfeld} Let $\g_{{\rm diag}} = \{(x, x) \mid x \in \g\}$. 
A Lagrangian splitting of the form $\gog = \g_{{\rm diag}} +\l$, where
$\l \in \Lgg$,
is called a Belavin-Drinfeld splitting. 
Let $G_{{\rm diag}} = \{(g, g) \mid g \in G\}$. It is shown in \cite{B-Dr}
(see also \cite[Corollary 3.18]{E-Lu2}) that every Belavin-Drinfeld
splitting of $\gog$ is conjugate by an element in $G_{{\rm diag}}$
to a splitting of the form 
\begin{equation}\lb{eq-BD}
\gog = \g_{{\rm diag}} + \l^{\prime\prime}_{S_2, T_2, d_2, V_2},
\end{equation}
where $(S_2, T_2, d_2)$ is a Belavin-Drinfeld triple in the sense that
\[
S_2^{d_2} = \{\alpha \in S_2 \mid d_2^n \alpha
\; \mbox{is defined and is in} \; S_2 \; \mbox{for} \; n = 1, 2, \cdots\}=\emptyset,
\]
and $V_2 \in {\mathcal L}_{{\rm space}}(\z_{S_2} \oplus \z_{T_2})$
is such that $\h_{{\rm diag}} \cap (V_2 + \{(x, \theta_{d_2}(x)) \mid x \in \h_{S_2}\}
= 0$. In other words, \eqref{eq-BD} is the special case of the splitting in 
\eqref{l1l2} with $\l_1 = \g_{{\rm diag}}$. Keeping the notation
as in \thref{rank-main}, we have $v_2 = 1$, and
the corank of $\Pi_{\l_1, \l_2}$ in $\O_1 \cap \O_2$ in this special case
simplifies to
\[
\dim \z_{S_2} + \dim \z_S 
-\dim(Y_2 \cap \Ad_{(\dot{w}_1, \dot{w}_2)} X_2)+
\dim (Z_2 \cap 
\Ad_{(\dot{w}_1, \dot{w}_2)} \widetilde{X}).
\]
When $\l^{\prime\prime}_{S_2, T_2, d_2, V_2} =\l_0: = \n^- \oplus \n + \h_{-{\rm diag}}$,
where $\n$ and $\n^-$ are respectively the span by positive and negative root
vectors and $\h_{-{\rm diag}} = \{(x, -x) \mid x \in \h\}$, the splitting
$\gog=\g_{{\rm diag}} + \l_0$ is called the {\it standard splitting} of
$\gog$ \cite{E-Lu2}. In this case, $N(\l_0) =B^- \times B$, where
$B^- = P_{\emptyset}^-$ and $B = P_\emptyset$
are two opposite Borel subgroups, and in the notation of \thref{rank-main},
$w_1 \in W$ and $w_2 \in W^T$. The corank of $\Pi_{\l_1, \l_2}$ in $\O_1 \cap \O_2$ 
in this special case further
simplifies to $\dim(\h_{-{\rm diag}} \cap \Ad_{(\dot{w}_1, \dot{w}_2)}\widetilde{X})$, where
\[
\widetilde{X} = \{(\Ad_{\dot{v}_1}^{-1}x, x) \mid x \in \z_{T(v_1)}, \; 
\theta_d \chi_S(\Ad_{\dot{v}_1}^{-1}x) = \chi_T(x)\} + \{(y, \theta_d(y))\mid
y \in \h_S\}.
\]
This formula 
has been obtained in \cite{E-Lu2}.
\eex

\subsection{The wonderful compactification of $G$}\lb{wond} 
Recall \cite{DP} 
that the wonderful compactification 
$\overline{G}$ of $G$ is the closure of the Lagrangian subalgebra 
$\g_{\mathrm{diag}}$ inside $\Lgg$. Let
$\g \oplus \g = \l_1 + \l_2$ be a Lagrangian splitting with
$\l_1$ and $\l_2$ given by \eqref{l1l2}. Then
$\overline{G}$ is a Poisson submanifold of $\Lgg$ 
with respect to the Poisson structure $\Pi_{\l_1, \l_2}$ because it is 
$(G \times G)$-stable. 
In \cite{LuY2}, to each of the above Lagrangian splittings we associated 
two partition ${\mathcal P}_i$, $i = 1, 2$, of $\overline{G}$ into 
finitely many 
smooth irreducible locally closed $N(\l_i)$-stable subsets. The strata
of ${\mathcal P}_i$ are
indexed by the Weyl group elements in \prref{N-orbits}
and are obtained by putting together the $N(\l_i)$-orbits 
corresponding to different continuous
parameters. When $\l_1 = \g_{{\rm diag}}$, the subsets in ${{\mathcal P}}_1$
are the $G_{{\rm diag}}$-stable pieces introduced by Lusztig 
\cite{L1, L2}. Each stratum of ${\mathcal P}_1$ and ${\mathcal P}_2$ is
a Poisson submanifolds of $(\overline{G}, \Pi_{\l_1, \l_2})$.
\thref{rank-main} shows that the corank of $\Pi_{\l_1, \l_2}$ 
at an $\l \in \overline{G}$ in $(N(\l_1)_\cdot \l)\cap (N(\l_2)_\cdot \l)$
depends only on the stratum of ${\mathcal P}_1$ (or ${\mathcal P}_2$)
to which $\l$ belongs. 

\end{document}